\documentclass[hidelinks,onefignum,onetabnum]{siamart220329}

\usepackage{lipsum}
\usepackage{amsfonts}
\usepackage{graphicx}
\usepackage{epstopdf}
\usepackage{algorithmic}
\ifpdf
  \DeclareGraphicsExtensions{.eps,.pdf,.png,.jpg}
\else
  \DeclareGraphicsExtensions{.eps}
\fi

\newsiamremark{remark}{Remark}
\newsiamremark{hypothesis}{Hypothesis}
\crefname{hypothesis}{Hypothesis}{Hypotheses}
\newsiamthm{claim}{Claim}

\headers{Solving HJB equations with FHT ansatz}{X. Tang, N. Sheng, L. Ying}

\title{Solving high-dimensional Hamilton-Jacobi-Bellman equations with functional hierarchical tensor}

\author{Xun Tang\thanks{First author. Institute for Computational and Mathematical Engineering, Stanford, CA 94305, USA. 
  (\email{xuntang@stanford.edu}).}
  \and
  Nan Sheng\thanks{Institute for Computational and Mathematical Engineering, Stanford, CA 94305, USA. 
  (\email{nansheng@stanford.edu}).}
  \and
  Lexing Ying\thanks{Department of Mathematics and Institute for Computational and Mathematical Engineering, Stanford, CA 94305, USA. 
  (\email{lexing@stanford.edu})}
  \funding{This work is partially supported by NSF grants DMS-2011699 and DMS-2208163. X.T. and L.Y. are supported by AFOSR MURI award FA9550-24-1-0254.}
  }

\usepackage{amsopn}

\makeatletter
\newcommand*{\addFileDependency}[1]{
  \typeout{(#1)}
  \@addtofilelist{#1}
  \IfFileExists{#1}{}{\typeout{No file #1.}}
}
\makeatother

\usepackage{amsmath,amsfonts,bm}

\def\eqref#1{equation~\ref{#1}}

\def\1{\bm{1}}

\DeclareMathAlphabet{\mathsfit}{\encodingdefault}{\sfdefault}{m}{sl}
\SetMathAlphabet{\mathsfit}{bold}{\encodingdefault}{\sfdefault}{bx}{n}

\newcommand{\R}{\mathbb{R}}

\DeclareMathOperator*{\argmin}{arg\,min}

\usepackage{amsmath,amssymb,a4wide,amsfonts}
\usepackage{latexsym}
\usepackage{indentfirst}
\usepackage{graphicx}
\usepackage{placeins}
\usepackage{enumerate}
\usepackage{booktabs}
\usepackage{algorithm}
\usepackage{algorithmic}
\usepackage{multirow}
\usepackage{optidef}
\usepackage{graphicx,color}
\usepackage{diagbox}
\usepackage[normalem]{ulem}

\ifpdf
\hypersetup{
  pdftitle={HJB-with-FHT},
  pdfauthor={X. Tang, N. Sheng, and L. Ying}
}
\fi

\newcommand{\edit}[1]{\textcolor{black}{#1}}

\begin{document}


\maketitle

\begin{abstract}
This work proposes a novel numerical scheme for solving the high-dimensional Hamilton-Jacobi-Bellman equation with a functional hierarchical tensor ansatz. We consider the setting of stochastic control, whereby one applies control to a particle under Brownian motion. 
In particular, the existence of diffusion presents a new challenge to conventional tensor network methods for deterministic optimal control.
To overcome the difficulty, we use a general regression-based formulation where the loss term is the Bellman consistency error combined with a Sobolev-type penalization term. We propose two novel sketching-based subroutines for obtaining the tensor-network approximation to the action-value functions and the value functions, which greatly accelerate the convergence for the subsequent regression phase. We apply the proposed approach successfully to two challenging control problems with Ginzburg-Landau potential in 1D and 2D with 64 variables.
\end{abstract}


\begin{keywords}
    High-dimensional Hamilton-Jacobi-Bellman equation; Semigroup method; Functional tensor network; Hierarchical tensor network; Curse of dimensionality.
\end{keywords}

\begin{MSCcodes}
    65M99, 15A69, 65F99
\end{MSCcodes}

\section{Introduction}

This paper focuses on solving the Hamilton-Jacobi-Bellman (HJB) equation:
\begin{equation}\label{eqn: HJB}
    \begin{cases}
        &\partial_{t}v(x, t) + \sup_{a \in A}\left[ b(x, a) \cdot  \nabla v(x, t)+ f(x, a)\right] + \frac{1}{\beta}\Delta v(x, t)  = 0,
    \quad x \in \mathcal{X} \subset \R^{d}, t \in [0, T),\\
    &v(y, T) = h(y), \quad  y \in \mathcal{X} \subset \R^{d},
    \end{cases}
\end{equation}
where $A\subset \R^{m}$ is the set of controls, $\beta$ is the inverse temperature controlling the diffusion coefficient, and the functions $b, f, h$ respectively encode the drift, running cost, and terminal cost of a stochastic control problem.

The HJB equation solution gives the optimal value function of its associated stochastic control problem. Let \((\Omega, \mathcal{F}, \left\{\mathcal{F}_t\right\}_{t \in [0, T]}, \mathbf{P})\) be a filtered probability space. The PDE in \Cref{eqn: HJB} models a stochastic control problem under the SDE dynamics:
\begin{equation}\label{eqn: langevin}
    dX_t = b(X_t, a_t) dt + \sqrt{2\beta^{-1}}dB_{t},
\end{equation}
where \(B_t\) is an \(\mathcal{F}_t\)-standard Brownian motion, and \(a_t \in A\) is an \(\mathcal{F}_t\)-adapted control. The cost functional \(J\) is defined by the following equation:
\begin{equation}
    J(t, y; \left(a_s \right)_{s \in [t, T]}) = \mathbb{E}\left[\int_{t}^{T}f(X_s, a_s) ds + h(X_T) \mid X_t = y\right],
\end{equation}
where the expectation is taken with respect to the trajectory of the Brownian motion.
The solution \(v\) to \Cref{eqn: HJB} can be represented as \(v(x, t) = \inf_{\left(a_s\right)_{s \in [t, T]}}\mathbb{E}\left[J(t, x, \left(a_s \right)_{s \in [t, T]})\right]\). In other words, \(v\) is the expected cumulative cost function under the optimal policy. We refer the readers to \cite{yong2012stochastic} for details. 

The stochastic control problem is a prevalent task in science and engineering with a wide range of applications. In this work, we focus on the general setting in which the drift term \(b\) is non-linear and the spatial variable \(x\) is in a high-dimensional space. In particular, this setting is challenging and precludes the application of traditional methods such as stochastic linear quadratic control \cite{chen1998stochastic}, finite volume methods \cite{wang2003numerical}, or Galerkin methods \cite{beard1997galerkin,bea1998successive}. 

In this work, we propose to use the parametric class of functional hierarchical tensor (FHT) \cite{hackbusch2012tensor,hackbusch2009new}, which has been recently applied to solving the Fokker-Planck equation \cite{tang2024solvinga} and Kolmogorov backward equations \cite{tang2025solving}. This parametric class of functions is suitable for reduced-order modeling under a hierarchical low-rank structure. In general, FHT has a good approximation guarantee under sufficient regularity characterized by a Sobolev-type norm \cite{schneider2014approximation}. This work employs the strategies used in deep reinforcement learning \cite{sutton2018reinforcement} and adapts them to function classes based on tensor networks. As the optimization subroutines differ drastically from neural network approaches, one of the main contributions of this work is the development of a full workflow for stochastic optimal control under a tensor network ansatz. 

Moreover, the use of the FHT function class is motivated by its applicability to model stochastic dynamics on function spaces, whereby the high-dimensional control problem comes from discretizing an infinite-dimensional Markov decision process. To give a motivating example, we consider the task of stabilizing a 2D field undergoing the gradient flow of a Ginzburg-Landau energy functional, a model widely used in the study of phase transitions and superconductivity. In the 2D G-L model, a ``particle" takes the form of a field \(x(z) \colon [0,1]^2 \to \R\). The Ginzburg-Landau potential is a functional $V(x)$ defined as follows:
\begin{equation}\label{eqn: GL potential infinite}
  V(x) = \frac{\lambda}{2}\int_{[0,1]^2} |\nabla_z x(z)|^2 \, dz + \frac{\mu}{4}\int_{[0,1]^2} |1- x(z)^2|^2 \, dz,
\end{equation}
where \(\lambda, \mu\) are parameters balancing the relative importance of the first and second term in \Cref{eqn: GL potential infinite}, which respectively represents the gradient energy and the double-well potential. One can apply a coarse control \(a \in \mathbb{R}\) with control direction \(\omega \colon [0, 1]^{2} \to \mathbb{R}\), and \(b\) admits the form
\begin{equation}\label{eqn: HJB drift example}
    b(x, a) = - \nabla V(x) + a\omega.
\end{equation}
Intuitively, the drift is in the direction of the gradient flow based on the Ginzburg-Landau functional, and the control applies a one-dimensional perturbation in the direction of \(\omega\), allowing for targeted influence on the system's evolution. 

Under discretization, one decomposes the unit square \([0,1]^2\) into a grid of \(d = m^2\) points $\{(ih,jh)\}$, for
$h=\frac{1}{m+1}$ and $1\le i,j\le m$. The discretized field is the $d$-dimensional vector $x=(x_{(i,j)})_{1\le i,j \le m}$, where \(x_{(i, j)} = x(ih, jh)\). The associated discretized potential energy is defined as
\begin{equation}\label{eqn: 2D GZ model}
V(x) = V(x_{(1,1)}, \ldots ,x_{(m,m)}) := \left[\frac{\lambda}{2} \sum_{w \sim w'}\left(\frac{x_{w} - x_{w'}}{h}\right)^2 +  \frac{\mu}{4} \sum_{w}\left(1 - x_w^2\right)^2\right]h^2, 
\end{equation}
where $w$ and $w'$ are Cartesian grid points and \(w \sim w'\) if and only if they are adjacent. Similarly, the control direction is discretized to the $d$-dimensional vector $\omega=(\omega_{(i,j)})_{1\le i,j \le m}$, and the drift term is again described by \Cref{eqn: HJB drift example}.

\subsection{Related work}\label{sec: related work}
\paragraph{Traditional methods in stochastic control}
Various classical methods have been developed to solve stochastic control problems. Examples include the finite volume method~\cite{richardson_numerical_2006,wang_numerical_2003}, the Galerkin method~\cite{beard1997galerkin,bea1998successive}, and the monotone approximation method~\cite{forsyth_numerical_2007}. For a detailed summary, we refer the readers to~\cite{rao_survey_2010} and reference therein. 
While these methods provide accurate numerical solutions, they cannot be applied to high-dimensional HJB equations due to the curse of dimensionality.

\paragraph{Neural network methods}
Great efforts have been made to solve high-dimensional HJB equations and related stochastic control problems using deep neural networks due to their effectiveness in high dimensions. One main approach is the deep backward stochastic differential equations (BSDE), where the BSDE formulation of HJB-type equations is parametrized and solved by deep neural networks \cite{e_deep_2017,han_solving_2018,nusken_solving_2021,pham_neural_2021,zhou_actor-critic_2021}.
Another approach is the physics-informed neural network~\cite{hu_tackling_2024,liu_pinn-based_2023}.
Other notable neural network approaches to tackle HJB equations include~\cite{han_deep_2016,nakamura-zimmerer_adaptive_2021,nakamura-zimmerer_causality-free_2020,darbon_neural_2021}.
While neural networks have been proven successful, 
algorithmic improvements might be difficult in general.

\paragraph{Tensor network methods}
Tensor network methods provide reduced-order parameterizations for high-dimensional functions and have seen great applicability in scientific computing. Several tensor-network-related approaches have appeared in the literature to solve HJB or related optimal control problems. Examples include the canonical polyadic decomposition~\cite{horowitz_linear_2014,kalise_polynomial_2018,chen_tensor_2022}, the tensor train decomposition~\cite{shetty_generalized_2023,dolgov_data-driven_2023,dolgov_tensor_2021,fackeldey_approximative_2022,gorodetsky_high-dimensional_2018,sommer_generative_2024,mosskull2020solving,eigel_dynamical_2023,gorodetsky_efficient_2015, oster_approximating_2021, oster2022approximating}.
In this work, we adopt the functional hierarchical tensor for optimal control problems. In contrast to most existing work where the problem is deterministic, this work focuses on addressing numerical challenges in a setting where control is applied to a particle under Brownian motion. In particular, the existence of diffusion precludes the common strategy of approximating the expected cost function through black-box interpolation, such as the approach taken in \cite{dolgov_data-driven_2023,shetty_generalized_2023}. \edit{Specifically, we choose the functional hierarchical tensor ansatz because it is well suited for representing value functions when the state variable corresponds to a discretized 2D field.}

\paragraph{Reinforcement learning}
Reinforcement learning (RL) has attracted significant attention for its great success in the past decade~\cite{sutton2018reinforcement}. 
The methodological development of RL techniques has laid a solid foundation
for solving optimal control problems. Notable RL methods include 
Q-learning~\cite{watkins1989learning,van2016deep} and the actor-critic method~\cite{konda1999actor,haarnoja2018soft}. 
Importantly, these methods have been extended
to solving stochastic optimal control problems under function approximation, examples of which can be found in~\cite{bachouch_deep_2022,wang_stochastic_2018,zhou_actor-critic_2021,han_actor-critic_2020}.
The empirical success of RL methods is partly based on solving the regression problem over the Bellman consistency error, and our work uses a similar loss function over functional tensor network architecture.

\subsection{Contents and notations} \label{sec: notation}

We outline the structure of the remainder of the manuscript. \Cref{sec: main workflow} goes through the algorithm for solving the HJB equation with the functional hierarchical tensor. \Cref{sec: FHT structure} introduces the functional hierarchical tensor structures used in this work. \Cref{sec: HJB time-stepping} details the optimization procedure used in the main workflow. \Cref{sec: numerical experiments} presents the numerical experiments and associated implementation details.

For notational compactness, we introduce several shorthand notations for simple derivations. For \(n \in \mathbb{N}\), let \([n] := \{1,\ldots, n\}\). For an index set \(S \subset [d]\), we let \(x_{S}\) stand for the subvector with entries from index set $S$. We use \(\Bar{S}\) to denote the set-theoretic complement of $S$, i.e. \(\Bar{S} = [d] - S\).

\section{Main formulation}\label{sec: main workflow}
\paragraph{Temporal discretization of the HJB equation}
To approximately solve the HJB equation in \Cref{eqn: HJB}, we propose to use a temporal discretization scheme whereby one selects time-steps \(0 = t_0 < t_1 < \ldots < t_K = T\) and restricts the control \(\left(a_s \right)_{s \in [0, T]}\) to be locally constant on the region \(s \in [t_k, t_{k+1})\). Moreover, this work uses the natural choice of uniform grid spacing with \(t_{k} = \frac{k}{K}T\). Under this discretization, the value function \(v(x, t)\) admits the interpretation of the optimal expected cumulative cost under locally constant policies for time \(s \in (t, T]\). Below $v(x,t_k)$ is often denoted as $v_k(x)$.

\paragraph{\(Q\) function}
Under the proposed discretization scheme, the central object used in our formulation is the \emph{action-value function}. We denote the action-value function for time \(t = t_{k}\) as \(Q_{k}\), defined as follows:
\begin{equation}\label{eqn: Q function}
    Q_k(x, a) = \mathbb{E}\left[\int_{s = t_{k}}^{t_{k+1}}f(X_s, a) ds + v(X_{t_{k+1}}, t_{k+1}) \mid \left( X_{t_k} = x \right)\wedge \left( \forall s \in [t_{k}, t_{k+1}),\,  a_{s} = a\right)\right].
\end{equation}
For the rest of this section, we explain the systematic procedure of using the \(Q_k\) function to obtain the value function \(v\). Subsequently, we specify the structure of \(Q_k\) and give an overview of how to obtain \(Q_k\).

\paragraph{Obtain \(v\) through backward time-stepping}
With temporal discretization, the goal of solving \Cref{eqn: HJB} is relaxed to finding the minimal expected cumulative cost function under locally constant controls. The proposed discretization scheme leads to a finite-horizon stochastic control problem with horizon \(K\). This work treats the stochastic control problem with the standard backward time-stepping strategy. For \(k = 0, \ldots, K-1\), assume that \(v(\cdot, t)\) has been solved for \(t \in \{t_{k+1}, \ldots, t_{K} = T\}\), and then the stochastic control problem at \(t = t_k\) can be solved with access to \(Q_{k}\) defined by \Cref{eqn: Q function}. Specifically, one can form the optimal control at \(t = t_k\) and state \(x\) by taking
\(a = \argmin_{\hat{a} \in A} Q_{k}(x, \hat{a})\), and likewise, one can calculate the value function at time \(t = t_{k}\) from \(  v_k(x)=v(x, t_{k}) = \min_{\hat{a} \in A} Q_{k}(x, \hat{a})\). 
Subsequently, one obtains \(Q_{k-1}\) for the \(t = t_{k-1}\), followed by solving \(v_{k-1}\) from \(Q_{k-1}\). Overall, one can continue the backward time-stepping strategy until \(t = t_{0}\) has been solved.

\paragraph{Obtain \(Q\) function through \(v\)}

To obtain the action-value functions for each \(t = t_k\), we propose to use the common \(Q\)-learning approach, for which we give a short introduction. One starts from selecting a collection of sampled state-action pairs \(\{x^{(i)}, a^{(i)}\}_{i = 1}^{N}\). For each state-action pair \((x, a)\), one evolves the SDE dynamics in \Cref{eqn: langevin} with a fixed control \(a\) for a duration of \(t_{k+1} - t_{k} = T/K\). After collecting the terminal state \(x'\) and the accumulated cost \(r\) for each state-action pair, one has the full collection of samples \(\{x^{(i)}, a^{(i)}, \left(x'\right)^{(i)}, r^{(i)}\}_{i = 1}^{N}\). Then, the \(Q\)-learning procedure for \(Q_k\) is obtained with a regression task over the Bellman consistency error:
\begin{equation}\label{eqn: Q learning}
    L(Q) = \frac{1}{N}\sum_{i = 1}^{N} \left\Vert Q(x^{(i)}, a^{(i)}) - v_{k+1}\left((x')^{(i)}\right) - r^{(i)} \right\Vert_{2}^{2},
\end{equation}
and \(v_k= v(\cdot, t_{k})\) can be obtained from \(Q_k\).

\paragraph{Function class of \(Q\) and \(v\)}
In this work, we propose to represent the action-value functions \(\{Q_k\}_{k=0}^{K-1}\) and value functions \(\{v_k\}_{k=0}^{K-1}\) through the functional hierarchical tensor ansatz. As the stochastic optimal control problem comes from the discretization of function spaces, the function space discretization leads to a lattice structure, which makes the FHT ansatz an attractive candidate for functional approximation. We defer the details of this ansatz to \Cref{sec: FHT structure}. 

\paragraph{Obtain continuous-in-time \(v\)}

After the series of \(K\) control problems have been solved following backward propagation, one obtains the collection of action-value functions \(\left\{Q_{k}\right\}_{k = 0}^{K-1}\) and value functions \(\left\{v_k := v(\cdot, t_{k}) \right\}_{k = 0}^{K-1}\). To obtain the continuous-in-time solution for \(v(\cdot, t)\) at non-grid points \(t\), one first finds the index \(k\) for which \(t \in [t_{k}, t_{k+1})\). Then, one can solve for the action-value function at time \(t\) following the definition in \Cref{eqn: Q function} with \(t\) in place of \(t_{k}\), and one then obtains \(v(\cdot, t)\) by taking minimization over the action-value function. A more efficient alternative is to directly take the linear interpolation between \(v_k\) and \(v_{k+1}\).

\paragraph{Summary of workflow}
We summarize the workflow here:
\begin{enumerate}
    \item Select \(K\) time steps \(0 = t_0 < t_1 < \ldots < t_K = T\) with \(t_{k} = \frac{k}{K}T\). 
    \item Run \(N\) independent stochastic simulations on the stochastic dynamic equation \Cref{eqn: langevin} on state-action pairs \(\{x^{(i)}, a^{(i)}\}_{i = 1}^{N}\) to time \(\frac{T}{K}\).
    \item Collect the associated terminal condition and the accumulated cost as \(\{(x')^{(i)}, r^{(i)}\}_{i = 1}^{N}\).
    \item Starting at \(k = K - 1\), obtain an FHT approximation of \(Q_{k}\) by regression over the Bellman consistency error in \Cref{eqn: Q learning} with the collected samples \(\{x^{(i)}, a^{(i)}, \left(x'\right)^{(i)}, r^{(i)}\}_{i = 1}^{N}\) and the previous time-step value function \(v_{k+1}\). Then, use \(Q_k\) to obtain \(v_k\). Iterate until one obtains \(Q_0, v_0\).
    \item For any \((x, t)\), output the continous-in-time solution either through computing the action-value function at time \(t\) or through linear interpolation.
\end{enumerate}




\section{Functional hierarchical tensor architecture}\label{sec: FHT structure}

The defining feature of this work is the functional approximation with the functional hierarchical tensor (FHT) architecture, and we give a short introduction of this ansatz following the more detailed exposition in \cite{tang2024solvinga}. As is seen in \Cref{sec: main workflow}, the FHT ansatz is used in the representation of three functions: (1) the value function \(v \colon \R^{d} \to \R\) which maps state to the expected cumulative cost, (2) the action-value function \(Q \colon \R^{d+m} \to \R\), which maps a state-action pair to the expected cumulative cost. Additionally, we consider (3) the Markov operator \(P \colon \R^{d+d+m} \to \R\), which models the probability transition following a given action. We shall give the corresponding hierarchical tensor structure for the three functions. 

\paragraph{Structure of the \(v\) function}
We first describe the functional tensor network representation of the \(d\)-dimensional value function \(v \colon \mathbb{R}^{d} \to \mathbb{R}\). 
Let \(\{\psi_{i;j}\}_{i = 1}^{n}\) denote a collection of orthonormal function basis over a single variable for the \(j\)-th co-ordinate, and let \(C_{v} \in \R^{n^{d}}\) be the tensor represented by a tensor network. The functional tensor network is the \(d\)-dimensional function defined by the following equation:
\begin{equation}\label{eqn: htn forward map}
    v(x)\equiv v(x_{1}, \ldots, x_{d}) = \sum_{i_{1}, \ldots, i_{d} = 0}^{n-1} \left(C_v\right)_{ i_1,\ldots, i_{d}} \psi_{i_1; 1}(x_1)\cdots \psi_{i_{d}; d}(x_{d}) = \left<C_v, \, \bigotimes_{j=1}^{d} \Vec{\Psi}_{j}(x_j) \right>,
\end{equation}
where \(\Vec{\Psi}_{j}(x_j) = \left[\psi_{1;j}(x_j),\ldots, \psi_{n;j}(x_j) \right]\) is an \(n\)-vector encoding the evaluation for the \(x_j\) variable over the entire \(j\)-th functional basis.
The functional hierarchical tensor structure of \(v\) is determined by the hierarchical tensor network ansatz of the tensor \(C_v\).

\begin{figure}[!ht]
    \centering
    \begin{minipage}{0.48\textwidth}
        \centering
        \includegraphics[width=\textwidth]{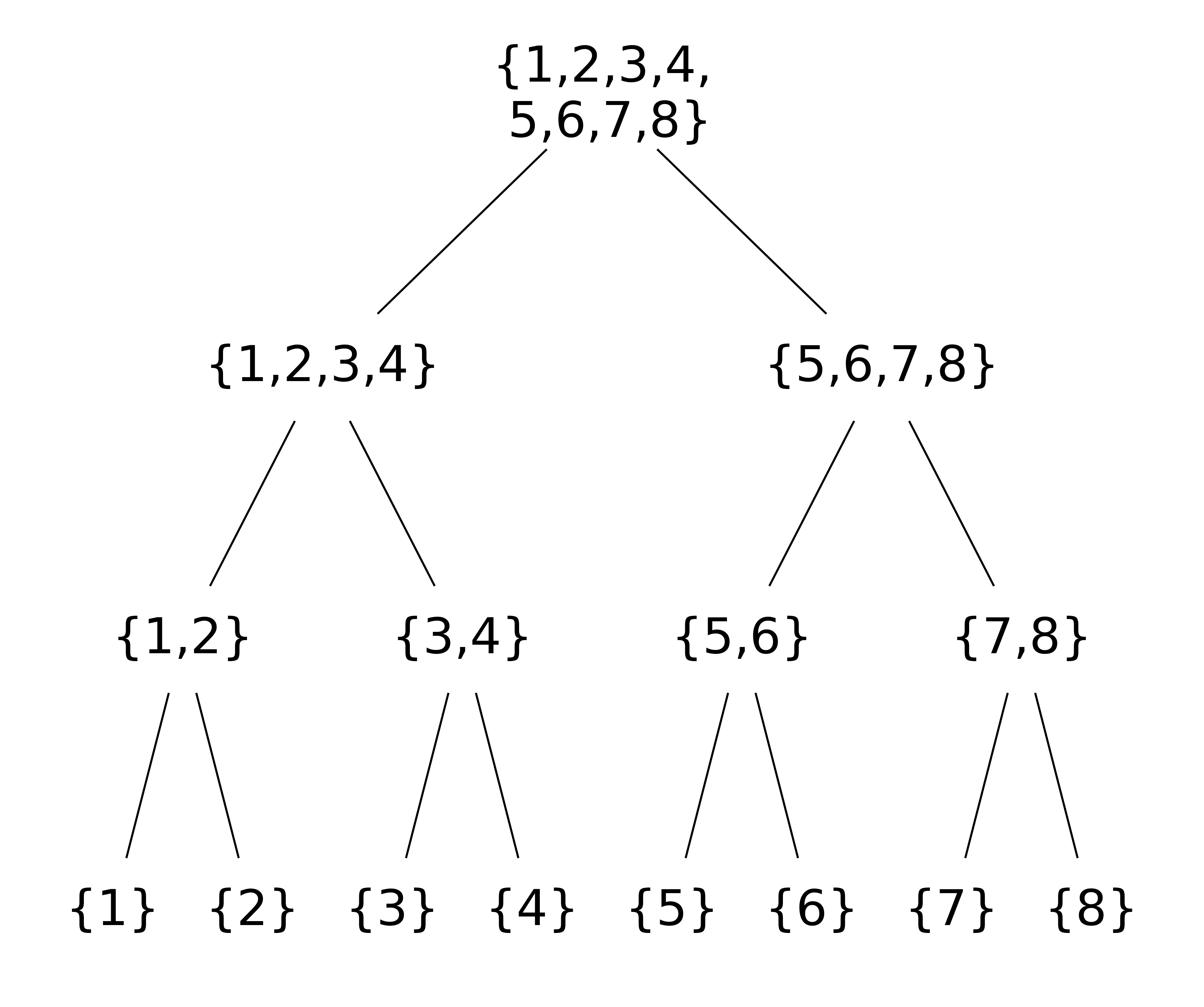}
        \caption{Binary decomposition of variables}
        \label{fig:binary_tree_8_nodes_subfig}
    \end{minipage}\hfill
    \begin{minipage}{0.50\textwidth}
        \centering
        \includegraphics[width=\textwidth]{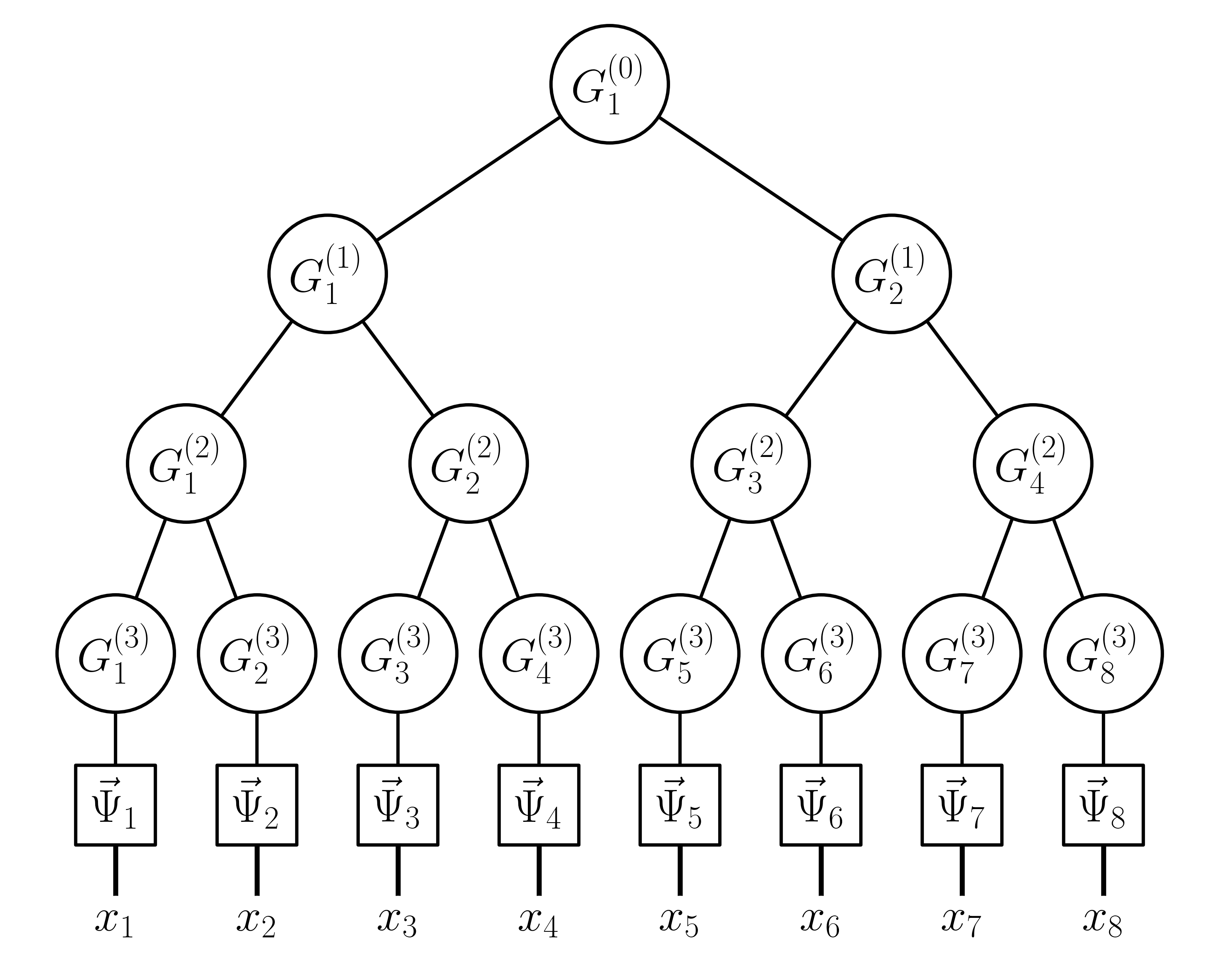}
        \caption{Functional hierarchical tensor structure of the value function}
        \label{fig:FHT_L_3}
    \end{minipage}
\end{figure}

A hierarchical tensor network is characterized by a hierarchical bipartition of the variable set. The exposition and notations largely follow from \cite{tang2024solvinga}. Without loss of generality, let \(d = 2^{L}\) so that the variable set admits exactly \(L\) levels of variable bipartition. As illustrated in \Cref{fig:binary_tree_8_nodes_subfig}, at the \(l\)-th level, the variable index set is partitioned according to 
\begin{equation}\label{eqn: bipartition}
    [d] = \bigcup_{k = 1}^{2^{l}} I_{k}^{(l)}, \quad I_{k}^{(l)} := \{ 2^{L - l + 1}(k-1) + 1, \ldots, 2^{L - l + 1}k\},
\end{equation}
which in particular implies the recursive relation that \(I_{k}^{(l)} = I_{2k-1}^{(l+1)}\cup I_{2k}^{(l+1)}\). 
Given the hierarchical bipartition, the tensor network structure of a hierarchical tensor network has the same tree structure as in the variable bipartition, which is illustrated by \Cref{fig:binary_tree_8_nodes_subfig}. In short, each internal bond of the tensor network is associated with an edge in the hierarchical binary decomposition, and each physical index is at the leaf node. 

\paragraph{Structure of the \(Q\) function}

\begin{figure}[!ht]
    \centering
    \begin{minipage}{0.48\textwidth}
        \centering
        \includegraphics[width=\textwidth]{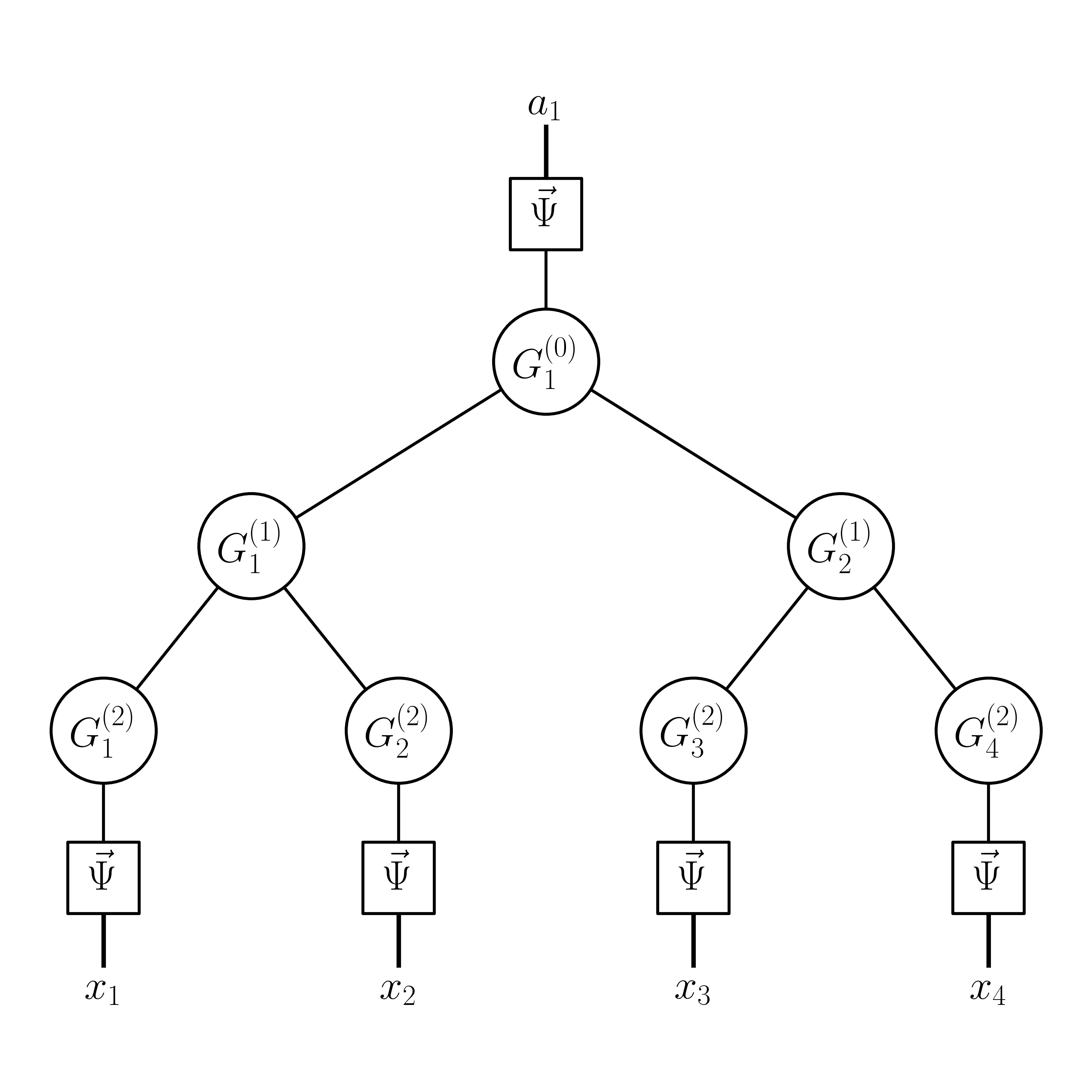}
        \caption{Functional hierarchical tensor structure of the action-value function for \(d = 4, m = 1\).}
        \label{fig:action_induced_FHT_L_3}
    \end{minipage}
    \begin{minipage}{0.48\textwidth}
        \centering
        \includegraphics[width=\textwidth]{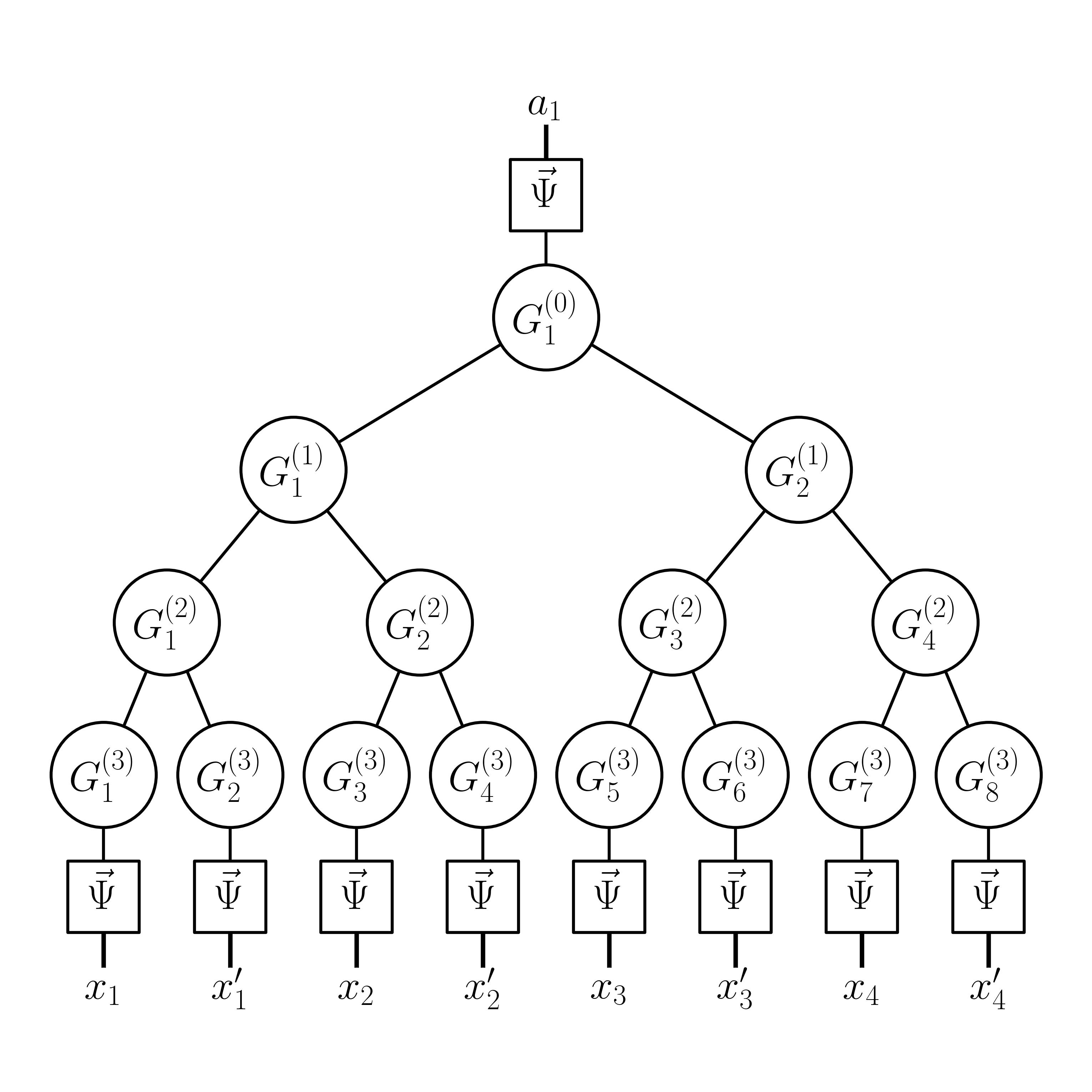}
        \caption{Functional hierarchical tensor structure of the Markov operator for \(d = 4, m = 1\).}
        \label{fig:action_induced_Markov_operator_L_3}
    \end{minipage}\hfill
\end{figure}

The input to the action-value function is the state-action pair, i.e., the concatenation of the \(d\)-dimensional state variable \((x_{1}, \ldots, x_{d})\) and the \(m\)-dimensional control variable \((a_1, \ldots, a_{m})\). Similar to the case of representing the value function, we propose to use a binary tree structure for the functional hierarchical tensor structure of the action-value function \(Q\). In this case, one needs to modify the placement of physical bonds to account for the added control variable. We propose to place physical bonds on the non-leaf nodes of the binary tree, one for each of the control variables. As tensor network models are universal approximations, non-optimal placement of variables can still enjoy good performance via suitably increasing the bond dimensions. In our numerical example, we consider the case of \(m = 1\), and therefore it is natural to place the physical bond at the root node of the binary tree, illustrated in \Cref{fig:action_induced_FHT_L_3}. \edit{We note that the proposed approach is heuristic, and the placement of control variables should be determined based on the specific structure of the problem at hand.}

\paragraph{Structure of the Markov operator}

The Markov operator \(P\) is defined such that \(P(x, x', a)\) defines the transition probability from state \(x\) to \(x'\) given one applies the action \(a\). The FHT structure of the Markov operator needs to be consistent with the FHT structure of the action-value function. The physical bonds for the \(x,x'\) variable are on the \((L+1)\)-th level of the binary tree ordered by the interlacing scheme \(z = (x_1, x'_1, \ldots, x_d, x'_{d})\). For the \(a\) variable, the placement of the physical bond is on the top \(L-1\) levels, and the location of each physical bond is the same as the respective location for the action-value function. In the case where \(m = 1\), one places the physical bond in the root node, illustrated in \Cref{fig:action_induced_Markov_operator_L_3}. We note that this network structure is similar to the Markov operator structure in \cite{tang2025solving}.

\section{Obtain value and action-value functions under the FHT ansatz}\label{sec: HJB time-stepping}

This section details the procedure for obtaining the FHT ansatz for the \(Q\) and \(v\) functions. Section \ref{sec: Markov operator} introduces the approach for \(Q\) function initialization with the Markov operator. Section \ref{sec: Q function ALS} goes through the approach for \(Q\) function regression. Section \ref{sec: value function ALS} explains how to obtain the \(v\) function. Section \ref{sec: form control} details how to use the \(Q\) function to form optimal control.

\subsection{Action-dependent Markov operator for \texorpdfstring{\(Q\)}{Q} function initialization}\label{sec: Markov operator}
One key component to the success of obtaining the action-value function approximation lies in a successful initialization. Here, we propose to use the Markov operator to obtain an initial value of the action-value function. \edit{The procedure employs a hierarchical sketching subroutine, as introduced in \cite{tang2025solving}.}

The expected cumulative cost for the state-action pair satisfies the following equation:
\[
Q_k(x, a) =\mathbb{E}\left[\int_{s = t_{k}}^{t_{k+1}}f(X_s, a) ds + v(X_{t_{k+1}}, t_{k+1}) \mid \left( X_{t_{k}} = x \right)\wedge \left( \forall s \in [t_{k}, t_{k+1}),\,  a_{s} = a\right)\right].
\]

To obtain an approximation of \(Q_k\), we define the \emph{Markov operator} \(P\) by the following equation:
\[
P(x, x', a) =\mathbb{P}\left[X_{t_{k+1}} = x' \mid \left( X_{t_k} = x \right)\wedge \left( \forall s \in [t_{k}, t_{k+1}),\,  a_{s} = a\right)\right],
\]
and so \(P\) is the transition probability from \(x\) to \(x'\) under the constant control \(a\). When \(\delta t = t_{k+1} - t_{k}\) is sufficiently small, one can take the approximation \(f(X_s, a) \approx f(X_{t}, a)\) in conjunction with the Bayes' rule to obtain the following equation:
\[
Q_k(x, a) \approx \int_{x'}P(x, x', a)v(x', t_{k+1}) dx' + f(x, a) \delta t, 
\] 
which only introduces an order \(\delta t^2\) error. Thus, one can see that modelling \(Q(x,a)\) can be performed by an integration formula between \(P\) and \(v(\cdot, t_{k+1})\). Similar to \cite{tang2025solving}, we propose to model \(P\) as a functional hierarchical tensor operator. Suppose that \(g, P\) has the structure as given in \Cref{sec: FHT structure} respectively, and then the integration between them can be done by efficient tensor contractions. By construction, the contraction of \(P, g\) has an FHT ansatz with the same structure proposed for the Q function. 

We now describe how the FHT-based approximation of \(P\) is obtained. Similar to the construction in \cite{tang2025solving}, we propose to approximate \(P\) with the Bayes rule. We let \(X, O\) be random state-action pairs following a chosen initial distribution \((X, O) \sim p_0\), and we let \(X'\) be the random variable formed by applying a constant control \(O\) from state \(X\) under the dynamics in \Cref{eqn: langevin} for a duration of \(\delta t\). One has
\begin{equation}\label{eqn: Markov operator by Bayes}
    P(x, x', a) = \frac{1}{p_0(x,a)}\mathbb{P}\left[X = x, X' = x', O = a \right],
\end{equation}
which shows that the conditional density \(P\) can be characterized through the probability density of the joint variable \((X, X', O)\). In our case, the state-action pair is supported on the bounded set \(\mathcal{X} \times A\), and it is natural to take \(p_{0} = \mathrm{Unif}(\mathcal{X}\times A)\). The uniform distribution leads to the formula
\begin{equation}\label{eqn: Markov operator new formula}
    P(x, x', a) = \mathrm{Vol}(\mathcal{X} \times A)\mathbb{P}\left[X = x, X' = x', O = a \right],
\end{equation}
and the joint density \(\mathbb{P}\left[X = x, X' = x', O = a \right]\) can be obtained by density estimation on the random variable \((X, X', O)\) under the FHT structure described in \Cref{sec: FHT structure}. As the approach is only a simple modification of the hierarchical sketching of the FHT ansatz in \cite{tang2025solving}, we include the implementation detail in \Cref{app: markov operator detail}.

\subsection{\texorpdfstring{\(Q\)}{Q} function regression under Bellman consistency error}\label{sec: Q function ALS}
After the action-value function \(Q_k\) is initialized, we can take the approach explained in \Cref{sec: main workflow} to refine the quality of action-value functions through alternating minimization. By selecting a collection of state-action pairs \(\{x^{(i)}, a^{(i)}\}_{i = 1}^{N}\), one evolves the dynamics in \Cref{eqn: langevin} and obtain training samples \(\{x^{(i)}, a^{(i)}, \left(x'\right)^{(i)}, r^{(i)}\}_{i = 1}^{N}\) for the downstream regression task.

Similar to the observations in \cite{oster2022approximating}, experiments show that using regression-based loss is not numerically quick to converge for practical problems. Therefore, we propose to use further smoothness regularization. In this case, the samples are normalized so that one has \(\{x^{(i)},a^{(i)}\} \in [-1, 1]^{d+1}\). For the regularization, we use the induced norm of the \((d+1)\)-way tensor product of the Sobolev space \( H^{1}([-1, 1]) \otimes \ldots \otimes H^{1}([-1, 1]) \), the details of which can be found in \Cref{app: ALS detail}. We denote the norm as \(\left\Vert  \cdot \right\Vert_{\mathrm{mix}}\). The Sobolev-regularized \(Q\)-learning procedure for \(Q_k\) has the following loss function:
\begin{equation}\label{eqn: Q learning final form}
    L_{\mu}(Q) = \frac{1}{N}\sum_{i = 1}^{N} \left\Vert Q(x^{(i)}, a^{(i)}) - v\left((x')^{(i)}, t_{k}\right) - r^{(i)} \right\Vert_{2}^{2} + \mu \left\Vert Q \right\Vert_{\mathrm{mix}}^2.
\end{equation}

\edit{We note that the mixed Sobolev norm \(\left\Vert \cdot \right\Vert_{\mathrm{mix}}^2\) is also used in \cite{oster2022approximating}, where it serves as a regularization term to promote smoothness in the resulting function.}
The regularization parameter \(\mu\) in \Cref{eqn: Q learning final form} is obtained by cross-validation, and a general rule is that \(\mu\) is such that the regularization term remains smaller than the regression loss term by orders of magnitude.

\paragraph{Optimization procedure}

In this work, we use the alternating least-squares (ALS) method \cite{holtz2012alternating} to optimize over the loss in \Cref{eqn: Q learning final form}. An important observation is that the mixed Sobolev norm term \(\left\Vert Q \right\Vert_{\mathrm{mix}}^2\) has a separable structure and admits an efficient tensor diagram. Therefore, the loss function in \Cref{eqn: Q learning final form} is quadratic in each tensor component. The ALS method takes advantage of the fact that the action-value function \(Q\) is parameterized in a multi-linear fashion, and the loss function \(L_{\mu}\) is a quadratic function when restricted to each tensor component. Therefore, the ALS method iteratively chooses tensor components and takes minimization over the loss function \(L_{\mu}\). 
Efficient implementation of the ALS scheme relies on reusing intermediate results on relevant tensor diagrams. By using a node schedule derived from the depth-first-search, one can ensure an efficient implementation where one optimizes over all of the tensor components with \(O(d)\) time. The complete implementation details are included in \Cref{app: ALS detail}.

\paragraph{Selection of training samples}
Similar to conventional supervised learning, the selection of the state-action pairs determines the quality of the action-value function. In this work, we use a simple sampling strategy whereby the samplings of states and actions are independent. The action \(a\) is sampled from a uniform distribution \(a \sim \mathrm{Unif}(A)\). For the state variable \(x\), half of the state samples \(x\) are sampled uniformly from \(\mathcal{X}\), and the other half is sampled from the dynamics in \Cref{eqn: langevin} with \(a = 0\). This sampling strategy ensures that there is training data for control-free trajectories. A more refined sampling strategy might lead to better performance.

\subsection{Obtaining the \texorpdfstring{\(v\)}{v} function}\label{sec: value function ALS}
After obtaining the action-value function \(Q_k\), obtaining the value function \(v_k\) falls within the category of function interpolation. For the action-value function \(Q_{k}\), we assume that the evaluation \(g(x) = \min_{a}Q_{k}(x, a)\) is efficient. As one can query the gradient \(\nabla_{a}Q_k(x, a)\) given the FHT ansatz of \(Q_k\), one can also use any conventional first-order optimization procedures (see \cite{nocedal1999numerical} for an overview). In our work, the numerical experiments consider the case in which \(a\) is one-dimensional, in which case any conventional optimization routine is sufficient to obtain the minimum. 

\paragraph{Initialization of the \(v\) function via hierarchical sketching}
To obtain an initialization of the value function, one important subroutine used in this work is functional hierarchical tensor interpolation. Through \(v_{k}(x) = \min_{\hat{a} \in A}Q_{k}(x, a)\), one can perform evaluation of \(v_k\). In essence, one performs repeated evaluations of \(v_k\) to obtain a sketched linear equation for the tensor component of the FHT structure of \(v_k\). The procedure is a novel use of the hierarchical sketching algorithm \cite{peng2023generative,tang2024solvinga}. Details are in \Cref{app: black-box interpolation detail}.

\paragraph{Alternating minimization for the \(v\) function}
After initialization, one selects a collection of state variables \(\{x^{(i)}\}_{i = 1}^{N}\), and one calculates the value function target by \(y^{(i)} = \min_{a}Q_{k}(x^{(i)}, a)\). The fine-tuned value function \(v_k\) is obtained by optimizing over the following Sobolev-regularized regression problem:
\begin{equation}
    L_{\mu}(v_k) = \frac{1}{N}\sum_{i = 1}^{N} \left\Vert v_k(x^{(i)}) - y^{(i)} \right\Vert_{2}^{2} + \mu \left\Vert v_k \right\Vert_{\mathrm{mix}}^2.
\end{equation}

Similar to the case of the action-value function, choosing the state variable sample determines the quality of the value function one obtains. For simplicity, this work uses the state variable samples from the state-action pair samples obtained in \Cref{sec: Q function ALS}.

\subsection{Generating stochastic optimal control with FHT ansatz}\label{sec: form control}
After one obtains the collection of action-value functions \(\left\{Q_{k}\right\}_{k = 0}^{K-1}\), one can use the forward propagation strategy to apply optimal control. At time \(t = t_{k}\), suppose that a particle is at state \(x\). With the strategies proposed in \Cref{sec: value function ALS}, one obtains the optimal control via \(a = \argmin_{\hat{a}}Q_{k}(x, \hat{a})\), and one can evolve the dynamics with the fixed control \(a\) for \(\delta t = t_{k+1} - t_{k}\). One then evolves the dynamics following the aforementioned strategy until one reaches \(t = T\).

\section{Numerical experiment}\label{sec: numerical experiments}
We apply the proposed method to a discretized Ginzburg-Landau model in 1D and 2D. For both experiments, the total time is \(T = 1\), and we take \(K = 10\) time steps with \(t_{k} = \frac{k}{10}\) for \(k = 0, \ldots, 10\). 

We describe the stochastic control problem considered for this case. For 1D and 2D Ginzburg-Landau models, the control-free (i.e. \(a = 0\)) stochastic dynamics has two meta-stable states \(y_{+} = (1, \ldots, 1)\) and \(y_{-} = (-1, \ldots, -1)\). We consider a stabilization task in which one applies control to the stochastic dynamics so that a particle is stabilized around the origin point \((0, \ldots, 0)\). Furthermore, we consider the coarse control case in which the control variable is one-dimensional with \(A = [-1, 1]\). 

To model this, we take the following cost functional:
\begin{equation}\label{eqn: cost functional GZ}
    J(t, y; \left(a_s \right)_{s \in [t, 1]}) = \mathbb{E}\left[\int_{s = t}^{1} \frac{1}{d}\left\Vert X_s \right\Vert^{2} + \left\Vert a_s\right\Vert^2  ds + \frac{1}{d}\left\Vert X_{1} \right\Vert^2 \mid X_t = y\right],
\end{equation}
which means that the instantaneous running cost \(f\) satisfies \(f(x, a) = \frac{1}{d}\left\Vert x \right\Vert^{2} + \left\Vert a \right\Vert^2\) and the terminal cost \(h\) satisfies \(h(y) = \frac{1}{d}\left\Vert y \right\Vert^{2}\). Thus, the model is set up so that the terminal cost and the instantaneous running cost depend on the particle's distance from the origin point \((0, \ldots, 0)\), and the instantaneous cost also has a quadratic term in the control variable \(a\). 

For the strength of diffusion, we take \(\beta = 1\), which makes the effect of diffusion significant without dominating the drift term coming from the Ginzburg-Landau functional.

\subsection{Implementation detail}
\paragraph{Initialization strategy}
In principle, the Markov operator initialization scheme for the action-value function and the FHT black-box interpolation algorithm for the value function can be used to initialize an FHT ansatz for any \(k = 0, \ldots, K-1\). In practical implementations, we use the proposed initialization scheme for initializing the action-value function and the value function only for \(k = K - 1\). Subsequent initialization is done by initializing the FHT ansatz \(Q_{k}, v_{k}\) with that of \(Q_{k+1}, v_{k+1}\). As \(t_{k+1} - t_{k}\) is small in our case, one has \(v_{k+1} \approx v_{k}\) and \(Q_{k+1} \approx Q_{k}\). By doing so, the cost for the more involved initialization procedure only needs to occur once for each ansatz, which simplifies the procedure while maintaining numerical accuracy. 

\paragraph{The choice of basis function}
Similar to \cite{tang2025solving}, the $\int_{z} |1 - x(z)^2|^2 dz$ term in the Ginzburg-Landau functional and the \(\left\Vert X_{s} \right\Vert^2\) term in the cost functional effectively ensures that the trajectory of a particle under an optimal trajectory is bounded within the domain \([-2, 2]^d\), which is why we take the domain to be the region \(\mathcal{X} = [-2, 2]^d\). Therefore, for the state variable \(x\), we use the polynomial representation by choosing a maximal degree parameter \(q\) and picking the first \(n = q + 1\) Legendre basis polynomial in $[-2,2]$. For the control variable \(a\), we similarly use the first \(n\) Legendre basis polynomial in $A = [-1, 1]$. In particular, we choose \(q = 6\). The effect of Sobolev-type regularization in \Cref{eqn: Q learning final form} implies that the presence of low-degree polynomials is more prevalent, which is why \(q = 6\) is sufficient in this case.

\paragraph{Bipartition structure}
The bipartition structure for the state variable exactly follows that of \cite{tang2024solvinga}, and we give an overview here for completeness. Assume for simplicity that a function \(x\) in $\Delta$-dimensions is discretized to the \(d=m^\Delta\) points denoted \(\{x_{i_1, \ldots, i_\Delta}\}_{i_1, \ldots, i_\Delta \in [m]}\) with \(m = 2^{\mu}\). For each index \(i_\delta\), one performs a length \(\mu\) binary expansion $i_\delta = a_{\delta 1}\ldots a_{\delta \mu}$. The variable \(x_{i_1, \ldots, i_\Delta}\) is given index \(k\) with the length \(\mu \Delta\) binary expansion \(k = a_{11} a_{21} \ldots a_{\Delta-1,\mu} a_{\Delta,\mu}\), and the associated binary decomposition structure follows from the given index according to \Cref{eqn: bipartition}. For other stochastic control problems, one can design a bipartition structure to fit the relationship between the variables.

\paragraph{Alternating minimization}
In the proposed alternating-least-squares (ALS) approach, one performs minimization over the tensor components once for each sweep schedule. One question is how many rounds of ALS are sufficient for the regression to be successful, as the Bellman consistency loss is a regression problem with noisy labels, and the model is under-parameterized. In general, one can use the typical regression strategy in supervised learning, where one performs training until the loss function does not improve. In our case, for every time step \(t_{k}\), we use five rounds of ALS for both the action-value function and the value function, which is sufficient for convergence. We remark that the success of ALS is heavily attributed to the initialization strategy proposed, and otherwise ALS cannot guarantee efficiency in training.

\subsection{1D Ginzburg-Landau model}

We first consider a 1D Ginzburg-Landau model. The potential energy is defined as
\begin{equation}\label{eqn: 1D GZ model}
V(x_1, \ldots ,x_m) := \left[\frac{\lambda}{2} \sum_{i=1}^{m+1}\left(\frac{x_{i} - x_{i - 1}}{h}\right)^2 + \frac{\mu}{4} \sum_{i = 1}^{m} \left(1 - x_i^2\right)^2\right]h,
\end{equation}
where \( h = \frac{1}{m+1} \) and \(x_0=x_{m+1}=0\). In particular, we fix \( m = 64 \), \( \lambda = 0.2\), \(\mu = 1 \). The dimension is \(d = 64\). For the control variable, we take a vector \(\omega \in \R^{d}\) where we let \(\omega_{i} = 1\) if \(\frac{i}{d} \in [0.25, 0.6]\), and \(\omega_{i} = 0\) otherwise. This is a discretization of the situation where one can apply a one-dimensional perturbation in the direction of \(\omega = \mathbf{1}_{[0.25, 0.6]}\). The drift term \(b\) admits the form
\begin{equation*}
    b(x, a) = - \nabla V(x) + 20a\omega,
\end{equation*}
and the cost functional \(J\) follows the form in \Cref{eqn: cost functional GZ}. 

\begin{figure}[!ht]
  \centering
\includegraphics[width = 0.6\textwidth]{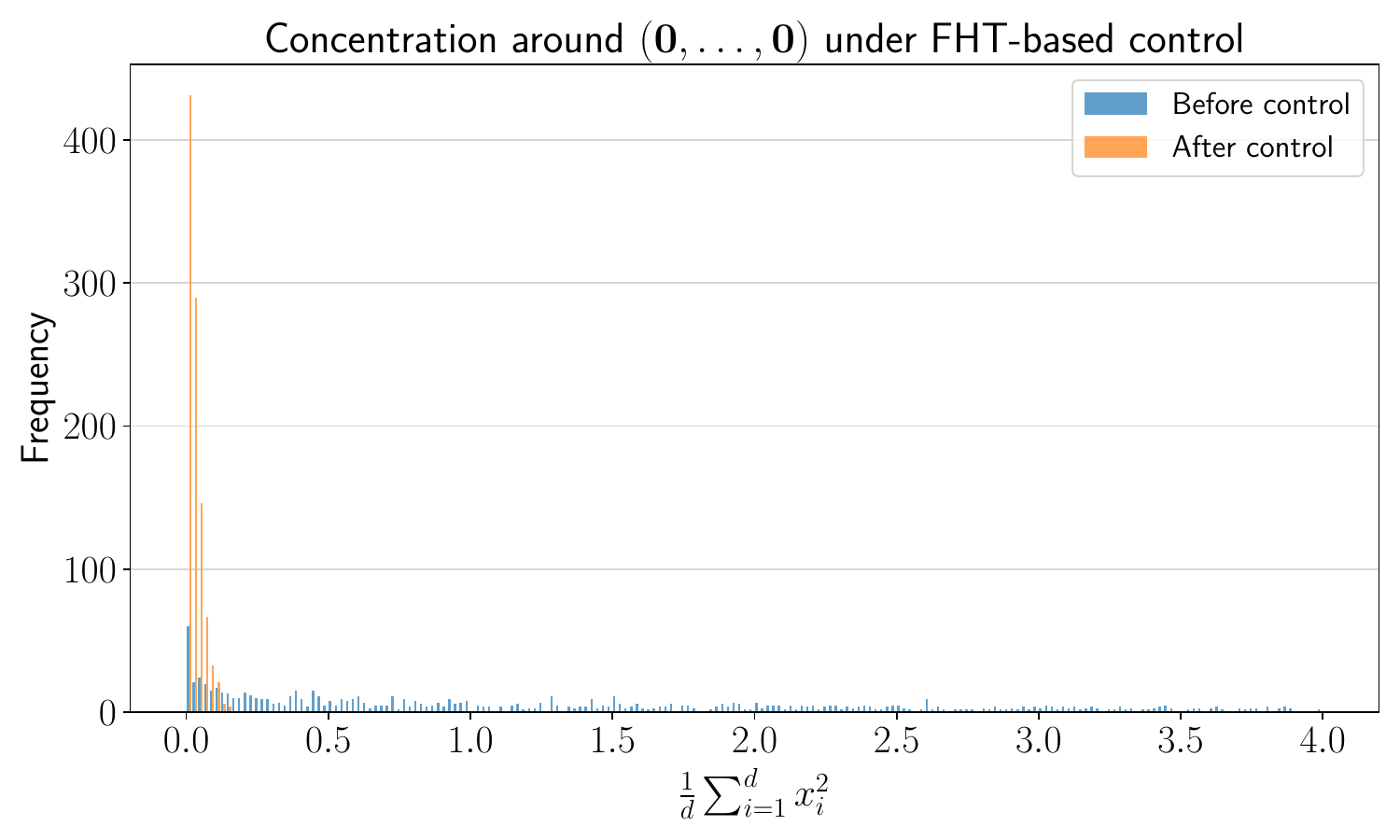}
  \caption{1D Ginzburg-Landau model. Distribution of the normalized distance of particles from the origin point. One can see that the proposed method succeeds at stabilizing a particle around the origin.}
  \label{Fig:1D_GZ_hist}
\end{figure}

In the first numerical test, we test whether a randomly generated particle can go to the origin instead of one of the meta-stable states \(y_{+}, y_{-}\). To test performance, we first generate \(500\) particles from the uniform distribution of \(\mathcal{X} = [-2, 2]^{d}\). In addition, we sample along another \(500\) points along the line segment \(I = \{(t, t, \ldots, t) \mid t \in (-2, 2)\}\) from the uniform distribution on \(I\). In \Cref{Fig:1D_GZ_hist}, we plot the concentration of the \(1000\) particles after applying locally constant stochastic control generated by the action-value functions. One can see that the FHT-based action-value function is successful at generating a control policy that guides the randomly generated particles to the origin.

\begin{figure}[!ht]
  \centering
\includegraphics[width = 1.0\textwidth]{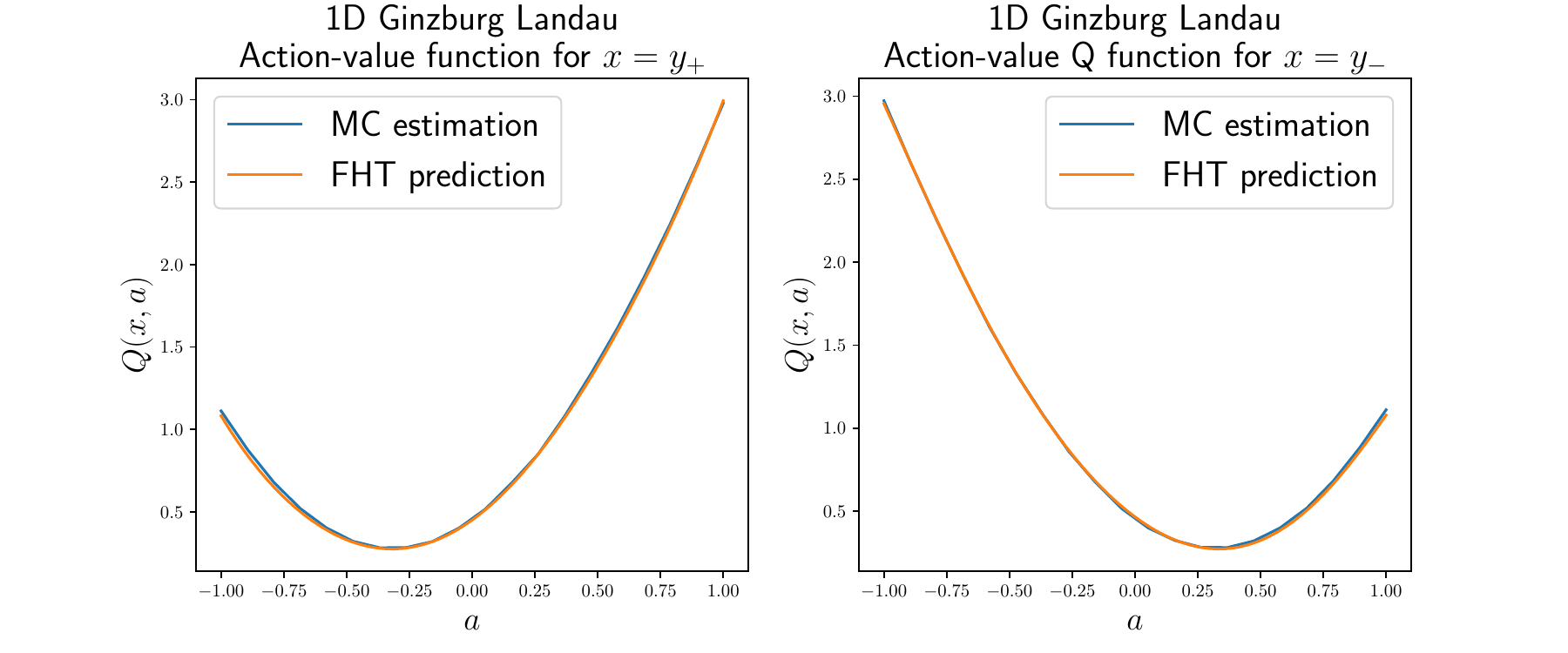}
  \caption{1D Ginzburg-Landau model. Plot of the action-value function at time \(t=t_{K-1}\). The ground truth action-value function is obtained through Monte-Carlo simulation.}
  \label{Fig:1D_GZ_Q_fun_plot}
\end{figure}

In the second numerical test, we test the accuracy of the action-value function the proposed method obtains. In particular, we consider the accuracy for the action-value function at time \(t = t_{K-1}\). This choice is because one can obtain the ground truth action-value function for any state-action pair by simply using Monte Carlo simulation. For illustration, we plot the action-value function \(Q_{K-1}(x, a)\) for varying values of \(a\) at two points \(x = y_{+}\) and \(x = y_{-}\). These two points are representative as the control-free trajectory stabilizes at the two meta-stable states, and it is desirable if the ansatz can form the optimal control from the two points accurately. In \Cref{Fig:1D_GZ_Q_fun_plot}, we plot the action-value function \(Q_{k}(x, a)\) for \(x = y_{+}\) and \(y_{-}\), which shows that the FHT approximation is very accurate. For both cases, the average relative error is \(1.6\%\). Importantly, one can see that the Sobolev regularization ensures that the action-value function is smooth across action values, which is appealing when one needs to use a first-order optimization approach to find \(\argmin_{\hat{a}}Q(x, a)\). The accuracy is similar when one chooses other values of \(x\), where we observe that the action-value function estimation remains accurate.

\begin{figure}[!ht]
  \centering
  \includegraphics[width = 0.5\textwidth]{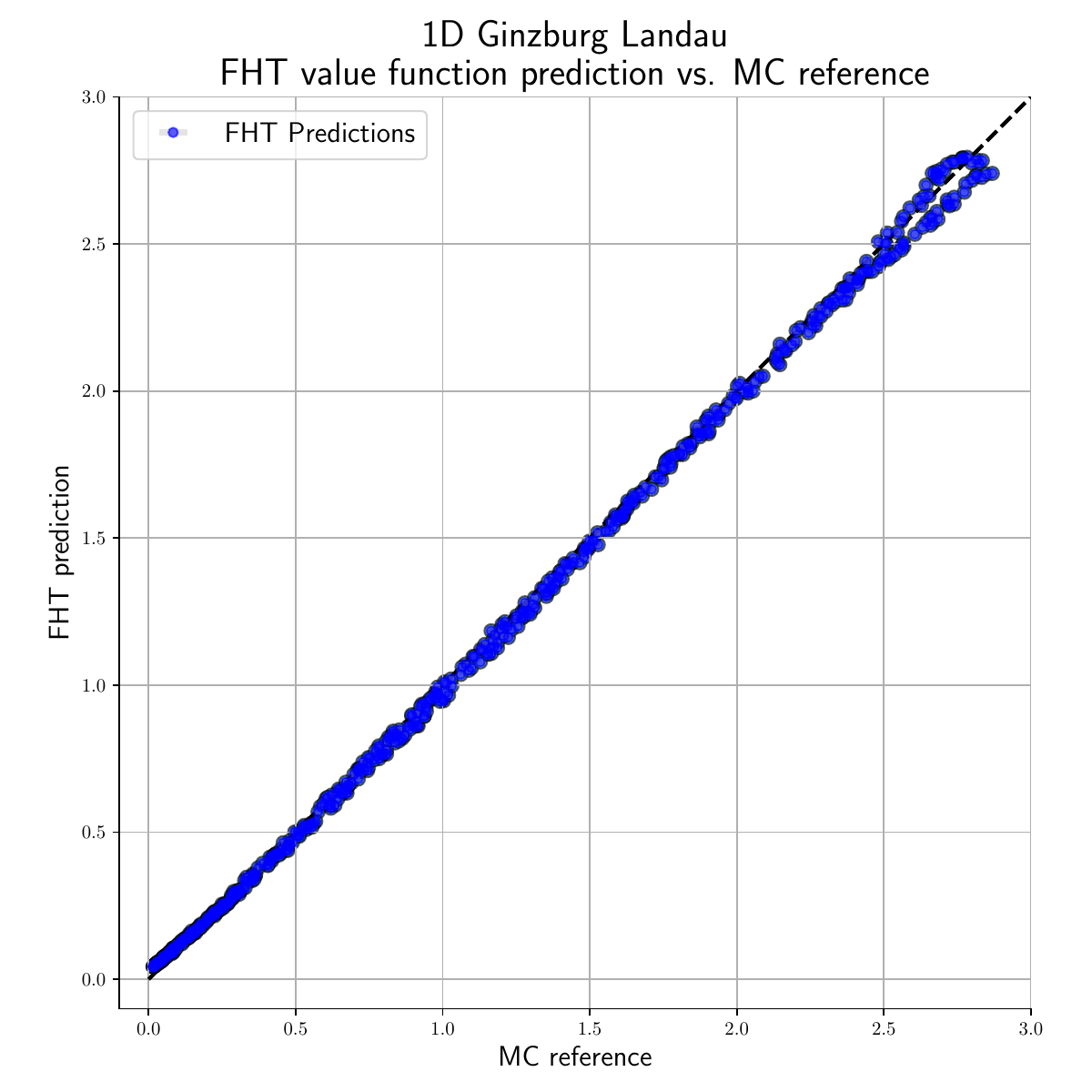}
  \caption{1D Ginzburg Landau model. Plots of the function evaluation of the FHT-based approximation of the value function at randomly sampled points. The reference values are obtained from Monte-Carlo estimations. }
  \label{Fig:1D_GZ_error_bar_plot_dist}
\end{figure}

In the third numerical test, we test the accuracy of the value function. Unlike the case for the action-value function, one does not have a ground-truth reference for the value function. Therefore, we similarly consider \(t = t_{K-1}\), and we propose to test the accuracy of \(v_{K-1}\) as the expected cumulative cost for the policy one forms using the action-value function at \(t = t_{K-1}\). In particular, we pick \(N = 600\) points, from which \(300\) points come from the uniform distribution of \(\mathcal{X} = [-2, 2]^{d}\), and the other \(300\) points come from the uniform distribution of the line segment \(I = \{(t, t, \ldots, t) \mid t \in (-2, 2)\}\). For each point, we sample \(100\) SDE trajectories following the optimal trajectories obtained from \(Q_{K-1}\), and we use the average cumulative cost to form the Monte-Carlo reference. As can be seen from \Cref{Fig:1D_GZ_error_bar_plot_dist}, the value function \(v_{K-1}\) matches well with references, and so the value function prediction is quite accurate. The results for other time steps are similar.

\subsection{2D Ginzburg-Landau model}

We then consider a 2D Ginzburg-Landau model. The potential energy follows from \Cref{eqn: 2D GZ model}.
In particular, we fix \( m = 8 \), \( \lambda = 0.5\), \(\mu = 1 \). The dimension is \(d = 64\). For the control variable, we take a vector \(\omega \in \R^{d}\) where we let \(\omega_{(i, j)} = 1\) if \((\frac{i}{d}, \frac{j}{d}) \in {[0.2, 0.8]}^2\), and \(\omega_{(i, j)} = 0\) otherwise. This is a discretization of the situation where one can apply a one-dimensional perturbation in the direction of \(\omega = \mathbf{1}_{{[0.2, 0.8]}^2}\). The drift term \(b\) admits the form
\begin{equation*}
    b(x, a) = - \nabla V(x) + 20a\omega,
\end{equation*}
and the cost functional \(J\) follows the form in \Cref{eqn: cost functional GZ}. 

\begin{figure}[!ht]
  \centering
\includegraphics[width = 0.6\textwidth]{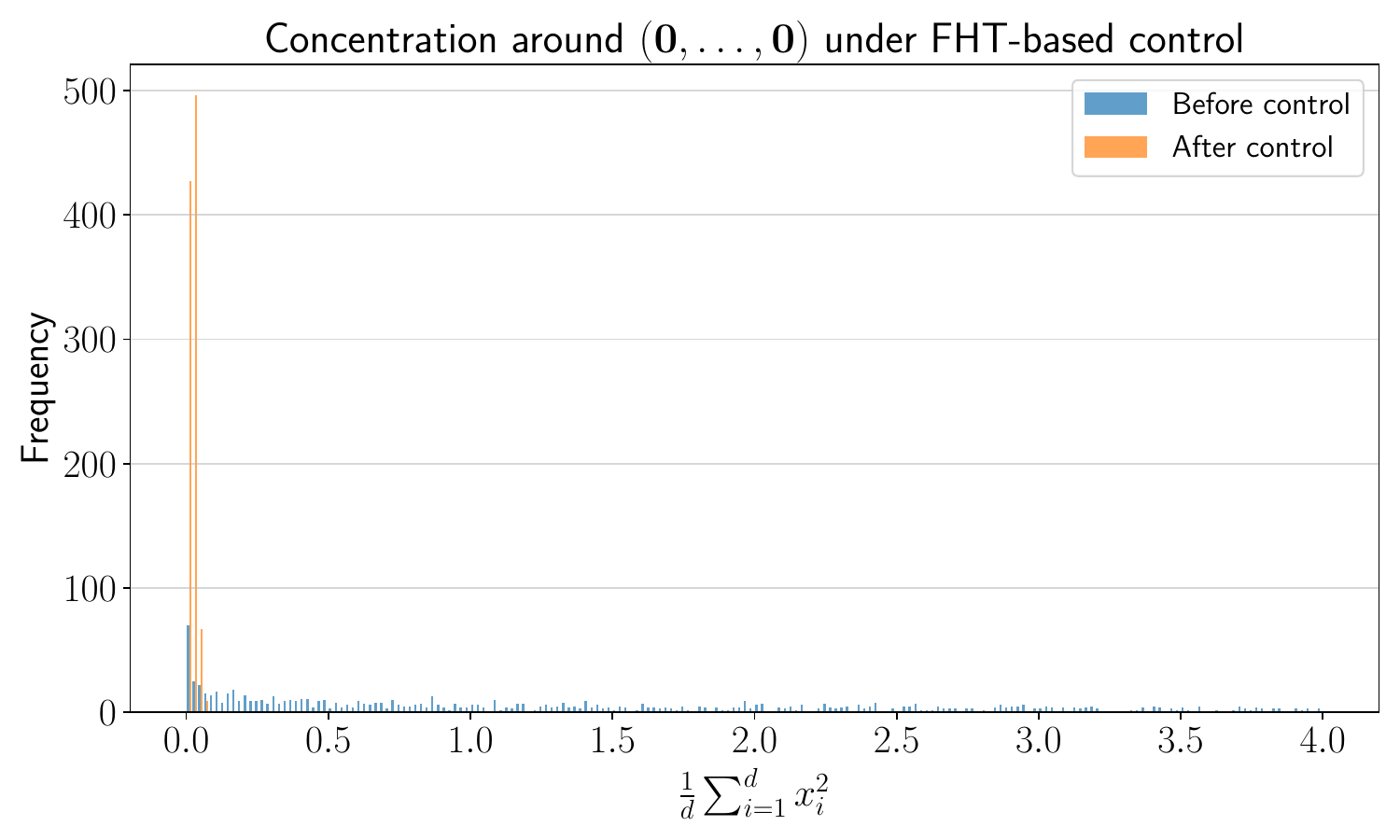}
  \caption{2D Ginzburg-Landau model. Distribution of the normalized distance of particles from the origin point. One can see that the proposed method succeeds at stabilizing a particle around the origin.}
  \label{Fig:2D_GZ_hist}
\end{figure}

\begin{figure}[!ht]
  \centering
\includegraphics[width = 1.0\textwidth]{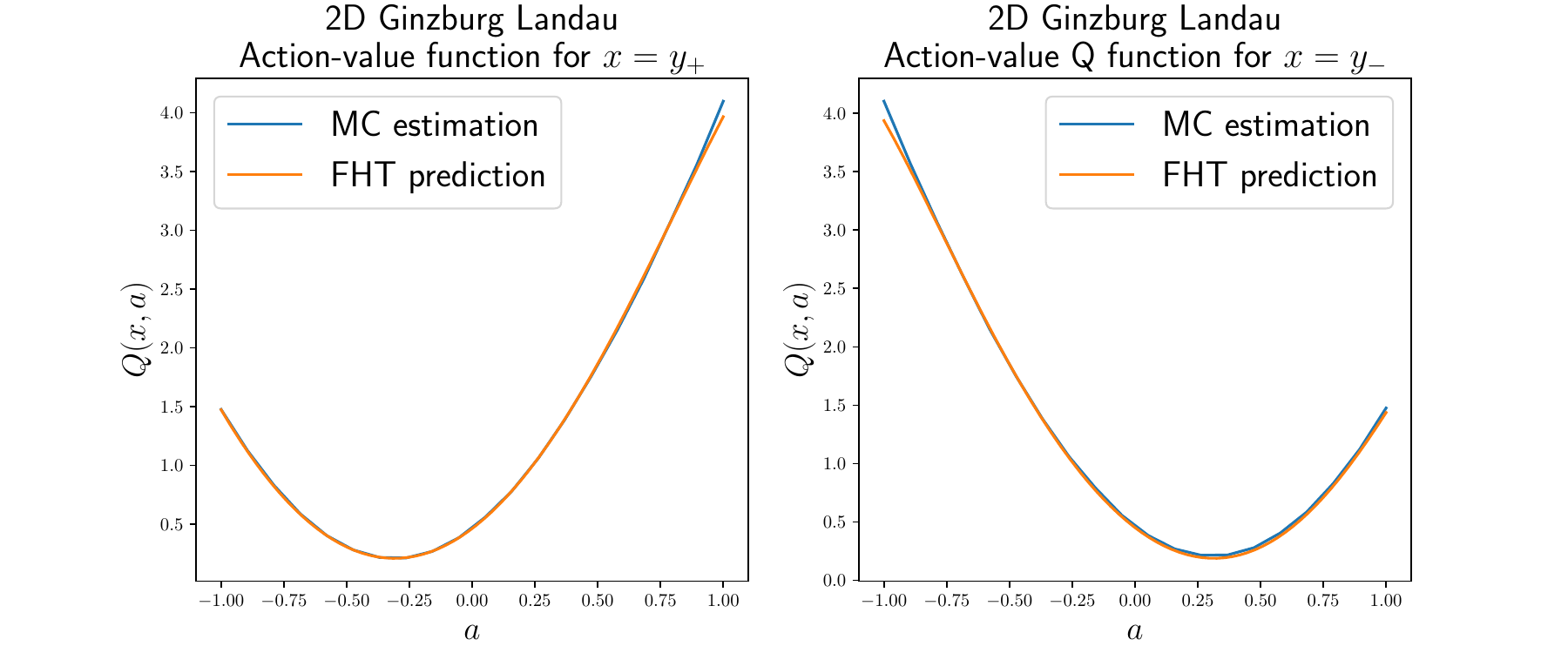}
  \caption{2D Ginzburg-Landau model. Plot of the action-value function at time \(t=t_{K-1}\). The ground truth action-value function is obtained through Monte-Carlo simulation.}
  \label{Fig:2D_GZ_Q_fun_plot}
\end{figure}

\begin{figure}[!ht]
  \centering
  \includegraphics[width = 0.5\textwidth]{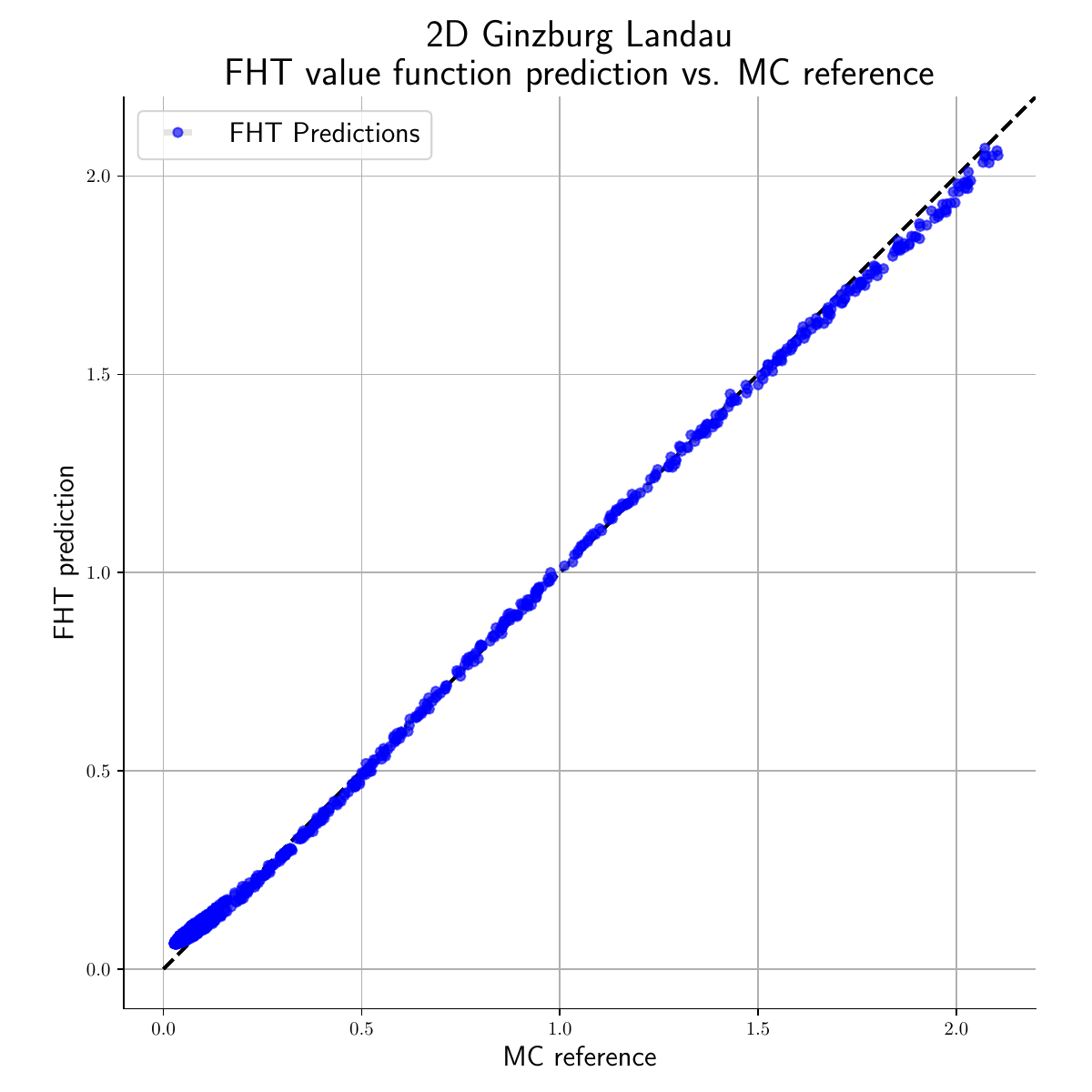}
  \caption{2D Ginzburg Landau model. The plot of the function evaluation of the FHT-based approximation of the value function at randomly sampled points. The reference values are obtained from Monte-Carlo estimations. }
  \label{Fig:2D_GZ_error_bar_plot_dist}
\end{figure}

We perform three numerical tests under the same setting as the 1D tests. 
In \Cref{Fig:2D_GZ_hist} we show the concentration of \(1000\) particles (\(500\) generated from the uniform distribution of \(\mathcal{X} = [-2, 2]^{d}\) and \(500\) along the line segment \(I = \{(t, t, \ldots, t) \mid t \in (-2, 2)\}\)) after applying locally constant stochastic control generated by the action-value functions. One can see that the FHT-based action-value function successfully takes the particles to the region around the origin. 
In \Cref{Fig:2D_GZ_Q_fun_plot}, we plot the action-value function at the two meta-stable states \(y_{+}, y_{-}\). The results show that the FHT approximate action-value function is close to the Monte-Carlo reference, and the average relative error in both cases is smaller than 3.4\%. The results for other state values are similar.
In \Cref{Fig:2D_GZ_error_bar_plot_dist}, we plot the value function compared to Monte Carlo references. The results show that the value function prediction similarly has good performance. 

\section{Conclusion}
We introduce a novel strategy to solve the Hamilton-Jacobi-Bellman equation with a functional hierarchical tensor network ansatz. The algorithm is applied to two problems with a discretized Ginzburg-Landau functional in 1D and 2D. This method demonstrates the effectiveness and high accuracy of using tensor-network-based methods to tackle high-dimensional PDEs. 


\bibliographystyle{siam_files/siamplain}
\bibliography{references}
\appendix
\section{Implementation details of alternating minimization}\label{app: ALS detail}

In this section, we detail the alternating minimization algorithm used in \Cref{sec: Q function ALS} and \Cref{sec: value function ALS}.  Without loss of generality, we discuss the case of \(g\) representing the value function, and one can extend to the action-value function case with mild change to the subroutine.

The goal is to find an approximate functional hierarchical tensor representation of a function $g(x)$ such that the following loss term is minimized:
\begin{equation} \label{eqn: als}
    L_{\mu}(g) = \frac{1}{N}\sum_{i = 1}^{N} \left\Vert g(x^{(i)}) - y^{(i)} \right\Vert_{2}^{2} + \mu \left\Vert g \right\Vert_{\mathrm{mix}}^2,
\end{equation}
where \(\{x^{(i)}, y^{(i)}\}_{i = 1}^{N}\) is the training data set for the Bellman consistency error. For the case of the \(Q\) function at \(t = t_{k}\), one takes collected samples of \(\{x^{(i)}, a^{(i)}, \left(x'\right)^{(i)}, r^{(i)}\}_{i = 1}^{N}\) by simulating the trajectory of \(x^{(i)}\) under action \(a^{(i)}\), and the regression target is \(y^{(i)} = r^{(i)} + v\left(\left(x'\right)^{(i)}, t_{k+1}\right)\). In the case of the \(v\) function at \(t = t_{k}\), one takes samples \(\{x^{(i)}\}_{i=1}^{N}\) and the regression target is \(y^{(i)} = \min_{\hat{a} \in A}Q_{k}(x^{(i)}, a)\).

The mixed Sobolev norm
\( \left\Vert \cdot \right\Vert_{\mathrm{mix}}\) is defined by  \[ \left\Vert f \right\Vert_{\mathrm{mix}} = \left( \sum_{k_1 =0}^{1} \cdots \sum_{k_d=0}^{1} \left\Vert \left(\frac{{\partial}}{{\partial}_{x_1}}\right)^{k_1} \cdots \left(\frac{{\partial}}{{\partial}_{x_d}}\right)^{k_d} f \right\Vert^2 \right)^{\frac{1}{2}}, \] 
where one can see that \( \left\Vert \cdot \right\Vert_{\mathrm{mix}}\) is the norm equipped for
the mixed Sobolev space \(H_{\mathrm{mix}}^{1} \equiv \otimes_{i=1}^d H^1 \). 
Under the functional hierarchical tensor representation, the function $g$ is parameterized by a collection of tensor components.
With all but one tensor component fixed, the loss function \(L_{\mu}\) in \Cref{eqn: als} is a quadratic function for a single tensor component. 
In essence, alternating minimization is a greedy approach aiming at solving each tensor node sequentially by minimizing the corresponding local quadratic function, with the goal of reaching the global minimum of \Cref{eqn: als}.

We first detail the derivation behind alternating minimization for a single tensor component.
We fix a node \(q\) that stands for the \(k\)-th block at level \(l\) of variable bipartition under the FHT ansatz. We consider \(G_{q}\), which stands for the tensor component corresponding to \(q\), and we show that \Cref{eqn: als} is a quadratic function in \(G_q\). For the rest of this section, we reserve the symbol \(G\) to denote the vectorization of \(G_q\), i.e., \(G := \mathrm{vec}(G_q)\).
For simplicity, we assume the general case $0< l < L$ so that \(q\) is neither the root nor the leaf node.

\paragraph{Bellman consistency error} 
First we show that the Bellman consistency error \(\left\Vert g(x^{(i)}) - y^{(i)} \right\Vert_{2}^{2}\) is a quadratic function of \(G_q\). The function $g$ in the FHT function class admits a representation
\begin{equation}\label{eqn: FHT fowrad pass}
    g(x)\equiv g(x_{1}, \ldots, x_{d}) = \left<C, \, \bigotimes_{j=1}^{d} \Vec{\Psi}_{j}(x_j) \right>,
\end{equation}
where \(\Vec{\Psi}_{j}\) is the basis set defined in \Cref{eqn: htn forward map} and $C$ is the coefficient tensor. The evaluation of \(g\) on the collected points \(\{x^{(i)}\}_{i = 1}^{N}\) reads \( g(x^{(i)}) = \left<C, \, \bigotimes_{j=1}^{d} \Vec{\Psi}_{j}(x_j^{(i)}) \right> \).

To show how the evaluation of \(g\) depends on \(G_q\), we let \(a := I_{2k-1}^{(l+1)}, b := I_{2k}^{(l+1)}\) and \(f = [d] - a \cup b\), where \(I_{k}^{(l)}\) follows hierarchical bipartition introduced in \Cref{eqn: bipartition}. The structural low-rankness property of the FHT ansatz implies that there exist \(r_a, r_b, r_f \in \mathbb{N}\), \(C_a \colon [n^{|a|}] \times [r_a] \to \R\), \(C_b \colon [n^{|b|}] \times [r_b] \to \R\), \(C_f \colon [r_f] \times [n^{|f|}] \to \R\) such that the following equations hold for \(G_q \colon [r_a] \times [r_b] \times [r_f] \to \R \):
\begin{equation} \label{eqn: uncrossed linear system for G_app_a}
    C(i_a,i_b,i_f) = \sum_{\alpha, \beta, \theta} 
    C_{a}(i_{a}, \alpha)C_{b}(i_{b}, \beta)G_{q}(\alpha, \beta, \theta)C_{f}(\theta, i_{f}),
\end{equation}
where terms such as \(i_a\) follows the multi-index notation as defined in \Cref{sec: notation}.

Combining \Cref{eqn: FHT fowrad pass} and \Cref{eqn: uncrossed linear system for G_app_a}, one has
\begin{equation} \label{eqn: A_aA_bA_f}
    g(x^{(i)}) = \sum_{\alpha, \beta, \theta}A_{a}(i, \alpha)A_b(i, \beta)A_f(i, \theta)G_{q}(\alpha, \beta, \theta),
\end{equation}
where the \(i\)-th row of matrices
\(A_a, A_b, A_f\) is the contraction of \(C_a, C_b, C_f\) with \(\bigotimes_{j \in a}\Vec{\Psi}_{j}(x_j^{(i)})\), \(\bigotimes_{j \in b}\Vec{\Psi}_{j}(x_j^{(i)})\) and \(\bigotimes_{j \in f}\Vec{\Psi}_{j}(x_j^{(i)})\). 
Now, we take the vectorization \(G = \mathrm{vec}(G_q)\) and let \( A \colon [N] \times [r_a \times r_b \times r_f] \to \R\) be the \emph{evaluation matrix} satisfying 
\begin{equation}
    A(i, (\alpha,\beta,\theta)) = A_{a}(i, \alpha)A_b(i, \beta)A_f(i, \theta),
\end{equation}
and then one has 
\(
\sum_{i = 1}^{n} \left\Vert g(x^{(i)}) - y^{(i)} \right\Vert_{2}^{2} = \left\Vert AG - y \right\Vert^2 \), which is indeed a quadratic function in the entries of \(G_q\).


\paragraph{Mixed Sobolev regularization}
Similarly, we show that the term \(\left\Vert g \right\Vert_{\mathrm{mix}}^2\) is a quadratic function in \(G_q\). When \(g\) admits an FHT parameterization in terms of \(C\) by \Cref{eqn: FHT fowrad pass}, one has \(\left\Vert g \right\Vert_{\mathrm{mix}}^2 = \mathrm{vec}(C)^{\top} K \mathrm{vec}(C) \) for an exponentially-sized matrix \(K\). An efficient evaluation of \(K\) relies primarily on the fact that \(K\) admits a separable tensor product structure. To see this, one can check the action of \(K\) on separable functions. For each site \(j \in [d]\), one can take an arbitrary function \(\phi_{j}\) within the vector space spanned by the basis function for the \(j\)-th variable \(x_j\). By taking \(f := \prod_{j = 1}^{d} \phi_{j}\), one sees that \(f\) is separable in the sense that \(f(x_1, \ldots, x_d) = \prod_{j=1}^{d} \phi_{j}(x_j)\) and
\[
\left\Vert \prod_{j = 1}^{d} \phi_{j} \right\Vert^2_{\mathrm{mix}} = \sum_{k_1 =0}^{1} \cdots \sum_{k_d=0}^{1} \left\Vert \left(\frac{{\partial}}{{\partial}_{x_1}}\right)^{k_1} \cdots \left(\frac{{\partial}}{{\partial}_{x_d}}\right)^{k_d}\prod_{j = 1}^{d} \phi_{j} \right\Vert^2 = \sum_{k_1 =0}^{1} \cdots \sum_{k_d=0}^{1} \prod_{j = 1}^{d} \left\Vert  \left(\frac{{d}}{{dx}}\right)^{k_j} \phi_{j}  \right\Vert^2,
\]
where the second equality holds due to Fubini's theorem. Subsequently, factorizing the exponentially sized summation allows us to rewrite
\[
\left\Vert \prod_{j = 1}^{d} \phi_{j} \right\Vert^2_{\mathrm{mix}} =  \prod_{j = 1}^{d} \left\Vert \sum_{k_j=0}^{1}  \left(\frac{{d}}{dx}\right)^{k_j} \phi_{j}  \right\Vert^2.
\]

Thus, the exponential-sized matrix \(K \colon [n_1 \times \ldots n_d ] \times [n_1 \times \ldots n_d ] \to \R\) can be written as a tensor product of matrices \(K = K_{1} \otimes \ldots \otimes K_d\), where \(K_{j} \colon [n_j] \times [n_j]\) is the matrix defined by
\[
K_{j}(i, i') = \left<\sum_{k = 0}^{1}\left(\frac{{d}}{{dx}}\right)^{k}\psi_{i; j}, \sum_{k' = 0}^{1}\left(\frac{{d}}{{dx}}\right)^{k'}\psi_{i'; j} \right>_{L^{2}([-1, 1])}.
\]

When restricted to the tensor component \(G_q\), one can likewise show that the Sobolev mixed norm term admits a form \( G^{\top} M G \). Let \(K_a = \bigotimes_{j \in a}K_{j} \), \(K_b = \bigotimes_{j \in b}K_{j} \), \(K_f = \bigotimes_{j \in f}K_{j}\). We form matrices \(M_{a}\), \(M_b\), \(M_f\) by \(M_{a} = C_{a}^{\top}K_{a}C_{a}, M_{b} = C_{b}^{\top}K_{b}C_{b}, M_{f} = C_{f}^{\top}K_{f}C_{f}\), where we utilize the matrix structure for \(C_a, C_b, C_f\) defined in this section. The term \(\left\Vert g \right\Vert_{\mathrm{mix}}\) can thus be written as \(\left\Vert g \right\Vert_{\mathrm{mix}} = G^{\top}MG\), where
\(M \colon [r_a \times r_b \times r_f] \times [r_a' \times r_b' \times r_f'] \to \R \) satisfies \(M = M_a \otimes M_b \otimes M_f\), i.e.
\begin{equation} \label{eqn: M_aM_bM_f}
    M((\alpha,\beta,\theta), (\alpha',\beta',\theta')) = M_{a}(\alpha, \alpha') M_{b}(\beta, \beta') M_f(\theta, \theta') .
\end{equation}

\paragraph{Single-node optimization in alternating minimization}
Combining with the results above, the quadratic problem for the tensor component $G_q$ is thus \[\min\limits_{G} \left\Vert AG - y \right\Vert_{2}^{2} + \mu G^T M G.
\]
The solution under the single-node optimization procedure is given by
\begin{equation} \label{eqn:G}
    G_q = \mathrm{Reshape}(G,r_a,r_b,r_f), \quad G = {(A^T A + \mu M)}^{-1} A^T y.
\end{equation}

\begin{figure}[h!]
    \centering
    \includegraphics[width = 1\textwidth]{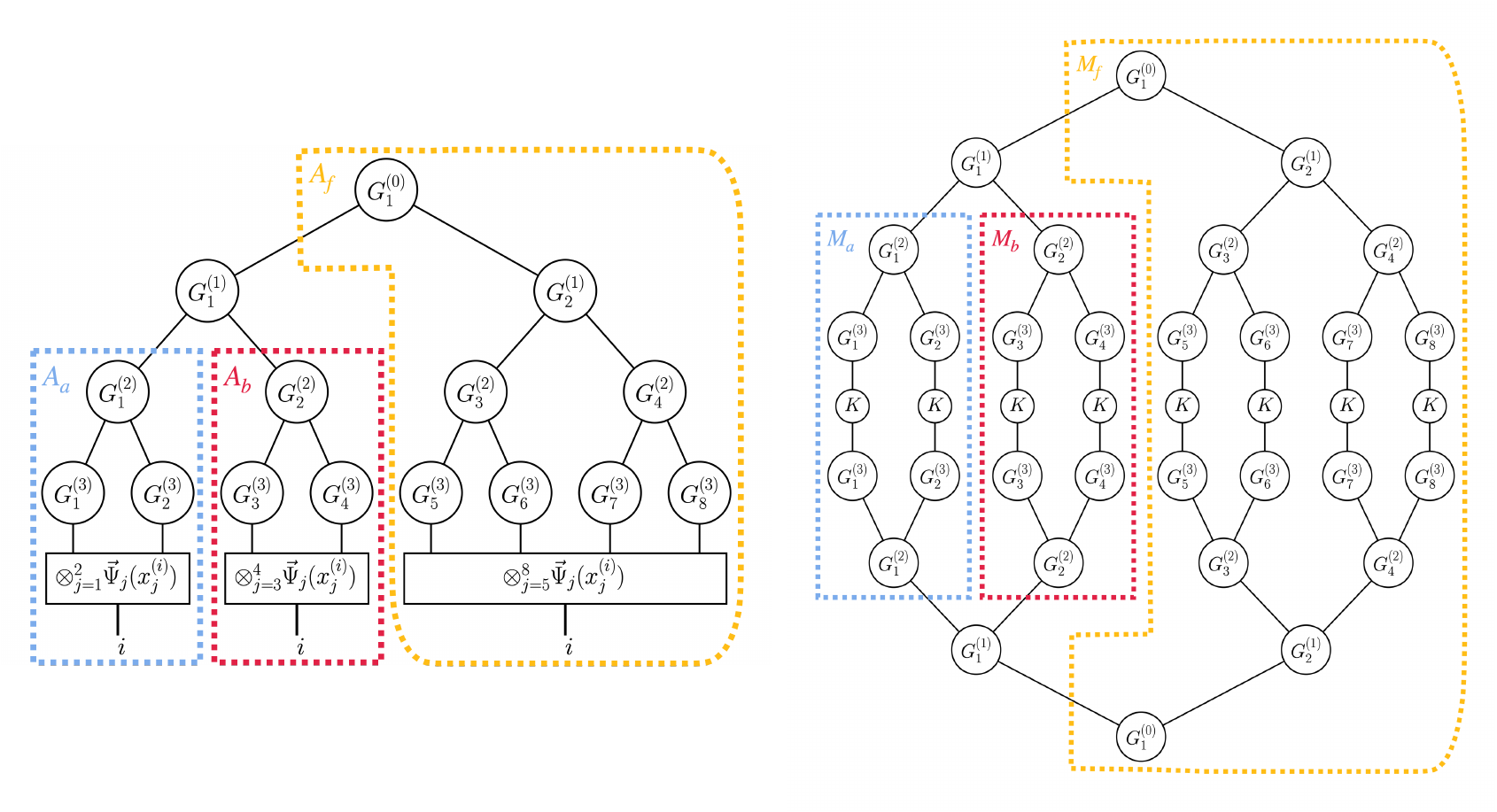}
    \caption{Left: diagrammatic expression of \(A_a\),\(A_b\),\(A_f\) appearing in \Cref{eqn: A_aA_bA_f} for \(G_1^{(1)}\). Right: diagrammatic expression of \(M_a\),\(M_b\),\(M_f\) appearing in \Cref{eqn: M_aM_bM_f} for \(G_1^{(1)}\).}
    \label{fig:ALS1}
\end{figure}

For an efficient numerical implementation, an important concept in performing single-node optimization for the FHT ansatz is the \emph{messages} sent to the node \(q\). The messages encapsulate the information necessary to compute the solution in \Cref{eqn:G}.
For illustration, we let \(q\) be the node corresponding to the tensor component \(G_1^{(1)}\) in \Cref{fig:ALS1} and explain the procedure.
To compute $G$ from \Cref{eqn:G}, we need to form tensors \(A_a,A_b,A_f\) and \(M_a,M_b,M_f\), and these six tensors are the messages to \(q\). 
For simplicity, we use \(m_{G_1^{(2)}\rightarrow G_1^{(1)}}\),\(m_{G_2^{(2)}\rightarrow G_1^{(1)}}\),\(m_{G_1^{(0)}\rightarrow G_1^{(1)}}\)
to denote
\((A_a,M_a)\), \((A_b,M_b)\), \((A_f,M_f)\) respectively in the example. We note that each edge of a functional hierarchical tensor corresponds to two messages in two directions, and there are \(O(d)\) messages in total.

\begin{figure}[h!]
    \centering
    \includegraphics[width = 0.5\textwidth]{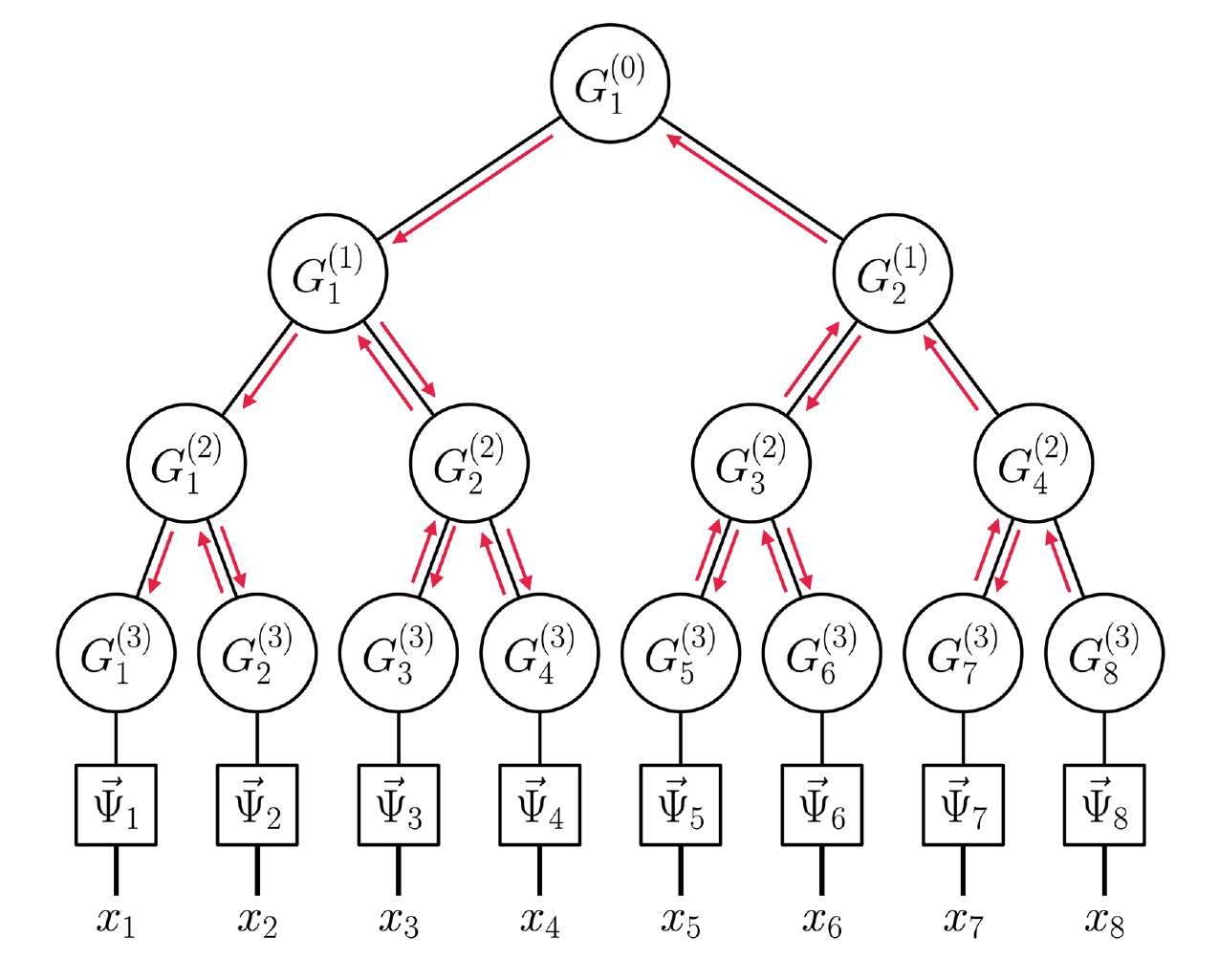}
    \caption{Depth-first-search path used in the sequential optimization of the functional hierarchical tensor.}
    \label{fig:ALS2}
\end{figure}

\begin{figure}[h!]
    \centering
    \includegraphics[width = 1\textwidth]{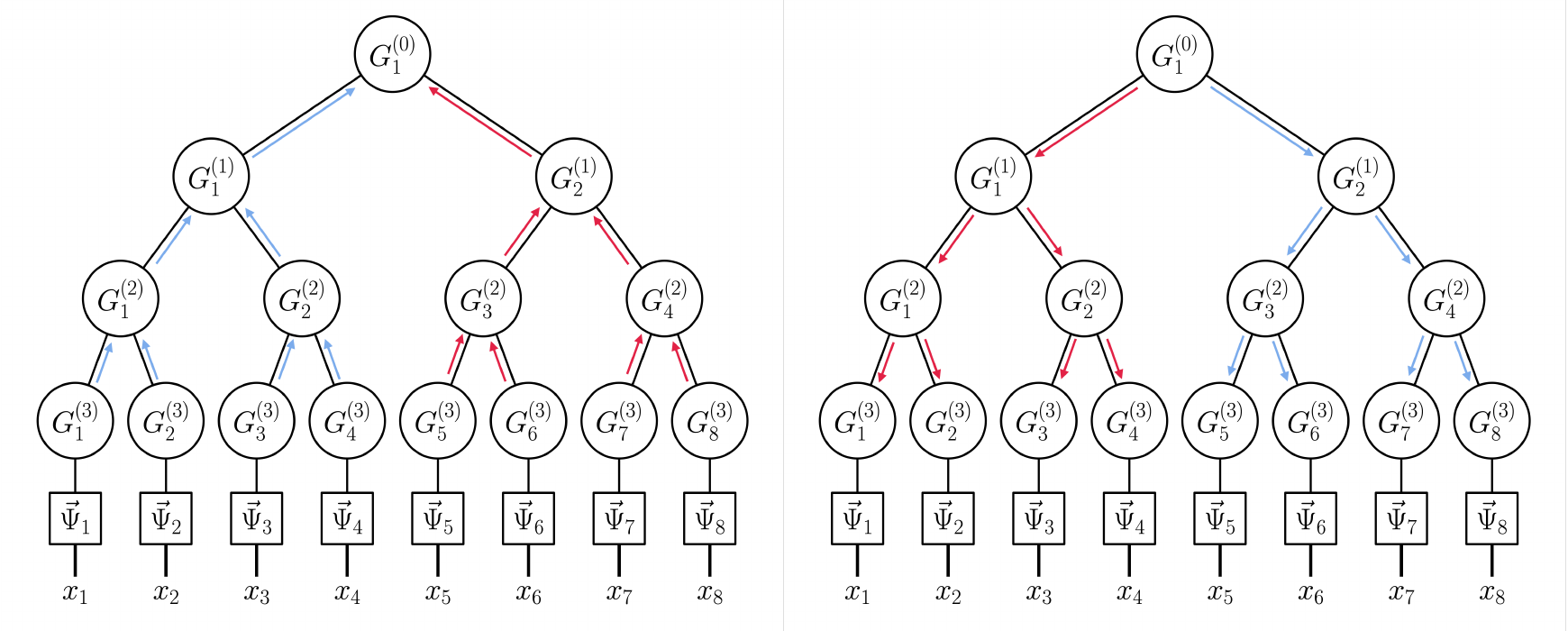}
    \caption{Message initialization process. Left: Leaf-to-root sweep process computing all children's messages. Right: Root-to-leaf sweep process computing all parent's messages.}
    \label{fig:ALS3}
\end{figure}

\paragraph{Node traversal for efficient alternating minimization}
We now go through the procedure for computing the optimal \(G_q\) under the alternating minimization procedure. Naively, one could compute the messages sent to each node \(q\) separately by sequential tensor contractions starting from the leaf nodes. It is not hard to see that this procedure would be inefficient, as it requires \(O(d)\) operations for updating a specific tensor component, leading to a total complexity of \(O(d^2)\) to update all of the tensor components.
In this work, we adopt a specific scheduling of tensor components optimization inspired by depth-first-search. The proposed strategy allows one to achieve an \(O(d)\) total complexity for node updates. The procedure is illustrated in \Cref{fig:ALS2} and \Cref{fig:ALS3}

Essentially, the proposed update strategy is based on the fact that one can update messages in \(O(1)\) cost for its nearest neighbor with a single tensor contraction. 
Using the example shown in \Cref{fig:ALS1}, one sees that \( m_{G_1^{(1)}\rightarrow G_1^{(2)}} \) can be computed with the knowledge of \( m_{G_1^{(0)}\rightarrow G_1^{(1)}} \), \( m_{G_2^{(2)}\rightarrow G_1^{(1)}} \) and \(G_1^{(1)}\). 
Specifically, \(A_f\) of \(G_1^{(2)}\) is given by a tensor contraction of three independent tensors: \(G_1^{(1)}\), \(A_b\) of \(G_1^{(1)}\) and \(A_f\) of \(G_1^{(1)}\); \(M_f\) of \(G_1^{(2)}\) is given by a tensor contraction of four independent tensors: \(G_1^{(1)}\), \(G_1^{(1)}\), \(M_b\) of \(G_1^{(1)}\) and \(M_f\) of \(G_1^{(1)}\).

The depth-first-search-inspired schedule we use is illustrated in \Cref{fig:ALS2} as a directed walk on the hierarchical tree.
We compute all of the messages in the beginning and start the walk from \(G_8^{(3)}\).
As three required messages are available, the tensor component \(G_8^{(3)}\) is optimized via \Cref{eqn:G}, 
and the schedule moves on to \( G_4^{(2)}\). 
To optimize \( G_4^{(2)}\), two required messages \(m_{G_7^{(3)}\rightarrow G_4^{(2)}}\) and \(m_{G_2^{(1)}\rightarrow G_4^{(2)}}\) are up to date. 
Since \(G_8^{(3)}\) has been updated, the message \(m_{G_8^{(3)}\rightarrow G_4^{(2)}}\) becomes out of date and one updates the message using the updated \(G_8^{(3)}\). Then \( G_4^{(2)}\) is optimized using the updated messages to it, and the schedule moves on to \(G_7^{(3)}\). To optimize \(G_7^{(3)}\), only \(m_{G_4^{(2)}\rightarrow G_7^{(3)}}\) is required, but \(m_{G_4^{(2)}\rightarrow G_7^{(3)}}\) has been out of date as both \(G_8^{(3)}\) and \(G_4^{(2)}\) have been updated. Since \(G_4^{(2)}\) has just been updated, the messages \(m_{G_2^{(1)}\rightarrow G_4^{(2)}}\) and \(m_{G_8^{(3)}\rightarrow G_4^{(2)}}\) are up to date and \(m_{G_4^{(2)}\rightarrow G_7^{(3)}}\) can be re-computed from tensor contractions of \(m_{G_2^{(1)}\rightarrow G_4^{(2)}}\), \(m_{G_8^{(3)}\rightarrow G_4^{(2)}}\) and \(G_4^{(2)}\). One then updates \(G_{7}^{(3)}\) using the updated messages to it. 

One continues the process and sequentially updates \(m_{G_{7}^{(3)} \rightarrow G_4^{(2)}}\) and \(m_{ G_4^{(2)} \rightarrow G_{2}^{(1)}}\). One sees that the messages to \(G_{2}^{(1)}\) are up to date, and one updates \(G_{2}^{(1)}\) and then moves on to update all descendants of \(G_{2}^{(1)}\) which are on the branch of \(G_{3}^{(2)}\). One can then update \(G_{1}^{(0)}\) and similarly update all descendants of \(G_{1}^{(0)}\) on the branch of \(G_{1}^{(1)}\), which then finishes the full update of nodes.

To realize the depth-first-search-inspired schedule,
all messages need to be initialized. We note that the cost of obtaining all messages is also \(O(d)\), and the total complexity of the algorithm is thus \(O(d)\). The approach of initializing all the messages comprised of two steps, whereby the first step is a leaf-to-root sweep computing all messages from a node to its parent, and the second step is a root-to-leaf sweep computing all messages from a node to its children. \Cref{fig:ALS3} illustrates our approach to initialize the messages. By sequential tensor contraction, in the leaf-to-root sweep, a node's message to its parent is computed using the messages from the node's children. In the root-to-leaf sweep schedule, one likewise uses sequential tensor contraction to compute the messages from a node to its children, which then allows one to compute the message from the children nodes to their own children nodes.

\begin{algorithm}[h]
\caption{Alternating minimization.}
\label{alg:1}
\begin{algorithmic}[1]
\REQUIRE Training data \(\{x^{(i)}, y^{(i)}\}_{i = 1}^{N}\).
\REQUIRE Chosen function basis \(\{\psi_{i;j}\}_{i = 1}^{n}\) for each \(j \in [d]\). 
\REQUIRE Regularization parameter \(\mu\).
\REQUIRE Number of iterations \(R\).
\STATE Gather all messages from tensor contractions using the two-step paths exemplified in \Cref{fig:ALS3}.
\FOR{iteration from 1 to \(R\)}
\FOR{each node \(q \in Q\) in the depth-first-search path exemplified in \Cref{fig:ALS2}}
    \STATE Update \(G_q\) using \Cref{eqn:G}.
    \STATE Update messages for the next node in sequence using a one-step tensor contraction.
\ENDFOR
\ENDFOR
\end{algorithmic}
\end{algorithm}


\paragraph{Algorithm summary}
The algorithm is summarized as \Cref{alg:1}, demonstrating the complete procedure for carrying out the alternating minimization algorithm for a given functional hierarchical tensor.  






\section{Implementation details of FHT black-box interpolation}\label{app: black-box interpolation detail}


We shall go over how to use the black-box interpolation method to obtain a functional hierarchical tensor representation of function $g(x)$ with access to querying its value at any given point $x$. For this subroutine, the input is a function \(g\) one has access to through evaluation, and the output is the FHT tensor representation of the function \(g\). \edit{We note that a similar idea has been explored in the tensor train setting~\cite{shi2024distributed, oseledets2010tt}.}

\begin{figure}[h!]
    \centering
    \includegraphics[width = 0.8\textwidth]{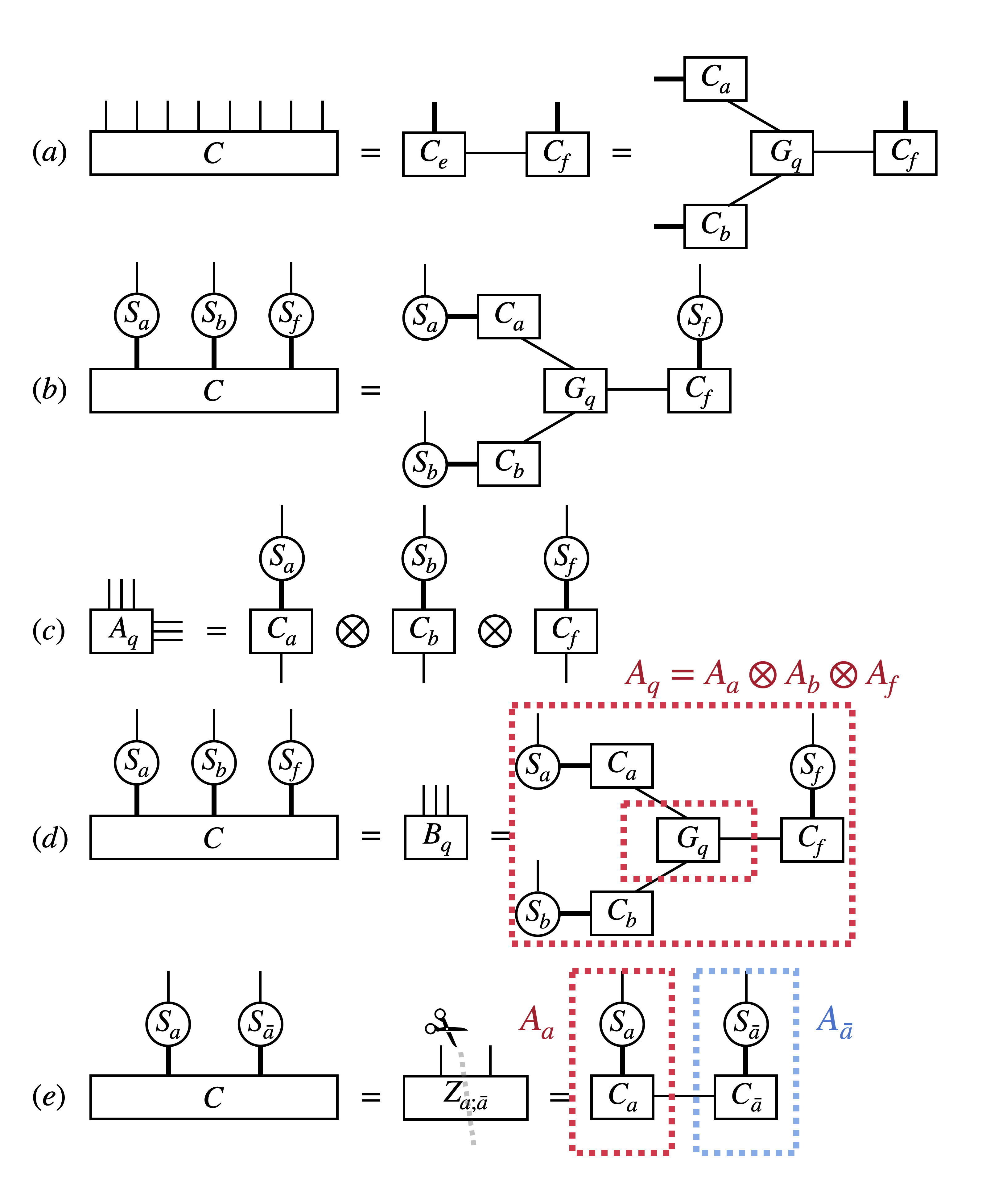}
    \caption{Tensor diagram representation of main equations used in \Cref{app: black-box interpolation detail}. \Cref{eqn: uncrossed linear system for G} is summarized in (a). \Cref{eqn: crossed linear system for G} is shown in (b). The coefficient term in \Cref{eqn: crossed linear system for G} is shown in (c). The definition of \(B_q\) is shown in (d). \Cref{eqn: eq for Z_ag} for \(Z_{a;\Bar{a}}\) is shown in (e), where the scissor symbol indicates the rank-\(r_{a}\) compression through SVD. }
    \label{fig:tensor_diagram_equation}
\end{figure}

The algorithm uses the assumption that $g(x)$ has a low-rank structure along a hierarchical bipartition, and it admits a functional hierarchical tensor representation (see Sec.~2 in Ref.~\cite{tang2024solvinga} for details). 
For the reader's convenience, important equations in this section are included in \Cref{fig:tensor_diagram_equation} in terms of the tensor diagram. Below, we discuss the main equation behind hierarchical black-box interpolation in the functional case. 


\paragraph{Equations} 

We use
the notations of the tensor diagram as shown in Fig. 3 for the functional hierarchical tensor. Again, we use $C$ to denote the tensor formed by the contraction of all tensor cores,
which satisfies 
\begin{equation}
    g(x)\equiv g(x_{1}, \ldots, x_{d}) = \left<C, \, \bigotimes_{j=1}^{d} \Vec{\Psi}_{j}(x_j) \right>.
\end{equation}
Here, we use similar notations to those defined in \Cref{eqn: htn forward map}.
The symbol
\(G_{q}\) denotes a single tensor core, where \(q\) is the node of the \(k\)-th block at level \(l\).
For simplicity, we assume the general case $0< l < L$ so that \(q\) is neither the root nor the leaf node. We define \(a := I_{2k-1}^{(l+1)}, b := I_{2k}^{(l+1)}\) and \(f = [d] - a \cup b\), where \(I_{k}^{(l)}\) follows notations introduced in \Cref{eqn: bipartition}. The structural low-rankness property of FHT implies that there exist \(r_a, r_b, r_f \in \mathbb{N}\), \(C_a \colon [n^{|a|}] \times [r_a] \to \R\), \(C_b \colon [n^{|b|}] \times [r_b] \to \R\), \(C_f \colon [r_f] \times [n^{|f|}] \to \R\) such that the following equations hold (illustrated in \Cref{fig:tensor_diagram_equation}(a)):
\begin{equation} \label{eqn: uncrossed linear system for G}
    C(i_a,i_b,i_f) = \sum_{\alpha, \beta, \theta} 
    C_{a}(i_{a}, \alpha)C_{b}(i_{b}, \beta)G_{q}(\alpha, \beta, \theta)C_{f}(\theta, i_{f}). 
\end{equation}

\Cref{eqn: uncrossed linear system for G} is exponential-sized, and one could not solve this equation directly. The hierarchical sketching algorithm~\cite{tang2024solvinga} introduces sketch functions to solve this over-determined linear system. In this algorithm, we sample structured random points and 
the evaluation of \(g\) on these points derives a sketched version of \Cref{eqn: uncrossed linear system for G}, which is the hierarchical sketching algorithm with a specific choice of sketch functions used in \cite{tang2024solvinga}. 

Let \(\tilde{r}_a, \tilde{r}_b, \tilde{r}_f\) be integers so that \(\tilde{r}_a > r_a, \tilde{r}_b > r_b, \tilde{r}_f > r_f\). To sketch out the exponential-sized tensors \(C_a, C_b, C_f\), we respectively select points \(\{x_{\mu}\}_{\mu \in [\tilde{r}_a]} \subset [-1, 1]^{|a|}\), \(\{x_{\nu}\}_{\nu \in [\tilde{r}_b]} \subset [-1, 1]^{|b|}\), \(\{x_{\zeta}\}_{\zeta \in [\tilde{r}_f]} \subset [-1, 1]^{|f|}\). The evaluation of \(g\) on the collected points leads to the following equation
\begin{equation}\label{eqn: FHT cross main equation}
    \forall \mu, \nu, \zeta, \quad g(x^{\mu}, x^{\nu}, x^{\zeta}) = \left<C, \, \bigotimes_{j \in a} \Vec{\Psi}_{j}(x^{\mu}_{j})  \bigotimes_{j \in b} \Vec{\Psi}_{j}(x^{\nu}_{j}) \bigotimes_{j \in f} \Vec{\Psi}_{j}(x^{\zeta}_{j}) \right>.
\end{equation}

To construct the linear system for \(G_{q}\), we first use the points \((x^{\mu}, x^{\nu}, x^{\zeta})\) to construct the right hand side term \(B_q\) of the linear system. In this case, the term \(B_q\) is obtained by repeated evaluation of \(g\) and given by the following equation:
\[
B_q(\mu, \nu, \zeta) := g(x^{\mu}, x^{\nu}, x^{\zeta}).
\]

We relate \(B_q\) to a sketched-down version of the exponential-sized \Cref{eqn: uncrossed linear system for G}. 
We introduce sketch tensors \(S_a \colon [n^{|a|}] \times [\tilde{r}_a] \to \R, S_b \colon [n^{|b|}] \times [\tilde{r}_b] \to \R, S_f \colon [n^{|f|}] \times [\tilde{r}_f] \to \R\) defined by 
\begin{equation}
    S_a(:,\mu) = \bigotimes_{j \in a}\Vec{\Psi}_{j}(x_j^{\mu}), \quad  S_b(:,\nu) = \bigotimes_{j\in b}\Vec{\Psi}_{j}(x_j^{\nu}), \quad S_f(:,\zeta) = \bigotimes_{j\in f}\Vec{\Psi}_{j}(x_j^{\zeta}),
\end{equation}
and one can readily see that \(B_q\) is obtained by contracting the \(C\) term in \Cref{eqn: uncrossed linear system for G} with \(S_a \otimes S_b \otimes S_f\).

Similarly, the left-hand-side of the linear system for \(G_q\) is naturally obtained by the contraction of \(C_a, C_b, C_f\) in \Cref{eqn: uncrossed linear system for G} with \(S_a, S_b, S_f\), which leads to the following equation:
\begin{equation}\label{eqn: crossed linear system for G}
    B_q(\mu, \nu, \zeta) = \sum_{\alpha, \beta, \theta}A_{a}(\mu, \alpha)A_b(\nu, \beta)A_f(\zeta, \theta)G_{q}(\alpha, \beta, \theta),
\end{equation}
which is illustrated in \Cref{fig:tensor_diagram_equation}(b) and \Cref{fig:tensor_diagram_equation}(d).
This linear system is thus
\[
(A_{a} \otimes A_{b} \otimes A_{f})G_{q} = B_{q}.
\]

The value of \(A_a, A_b, A_f\) depends on the gauge degree of freedom in \(C_a, C_b, C_f\), and the detailed procedure of systematically fixing the gauge is introduced in \cite{peng2023generative,tang2024solvinga}. Here, we briefly explain the procedure. Essentially, we show that one can obtain \(A_a\) without having access to \(C_a\). For \(\Bar{a} := [d] - a\), the structural low-rankness of \(g\) implies the existence of \(C_{\Bar{a}}\) such that there exists a linear system \(C(i_1, \ldots, i_d) = \sum_{\alpha}C_a(i_a, \alpha)C_{\Bar{a}}(\alpha, i_{{\Bar{a}}})\).
One can likewise obtain a sketched-down linear system by repeated evaluation of \(g\). 

Let \(\tilde{r}_{\Bar{a}}\) be an integer such that \(\tilde{r}_{\Bar{a}} > r_a\). 
Denoting $\mu,\phi \in [\tilde{r}_a],[\tilde{r}_{\Bar{a}}]$, we likewise define sketch points as \((x^\mu,x^\phi)\) and sketch tensor \(S_{\Bar{a}} \colon [n^{|a|}] \times [\tilde{r}_a] \to \R\) as 
\(S_{\Bar{a}}(:,\phi) = \bigotimes_{j \notin a}\Vec{\Psi}_{j}(x_j^{\phi})\). The function $g$ evaluated at the collection of points is given by
\begin{equation}\label{eqn: eq for Z_ag}
    Z_{a;{\Bar{a}}}(\mu, \phi) := g(x^{\mu}, x^{\phi}),
\end{equation}
and subsequently one has (illustrated in \Cref{fig:tensor_diagram_equation}(e))
\[
Z_{a;{\Bar{a}}}(\mu, \phi) = \sum_{\alpha}A_a(\mu, \alpha)A_{\Bar{a}}(\alpha, \phi)
\]
where \(A_{{\Bar{a}}}\) is the contraction of \(C_{\Bar{a}}\) with \(S_{\Bar{a}}\). 
Then, due to a gauge degree of freedom in choosing \(C_a\), there also exists a gauge degree of freedom in choosing \(A_a\). Thus, one can perform singular value decomposition (SVD) on \(Z_{a;{\Bar{a}}}\), obtaining \(Z_{a;{\Bar{a}}}(\mu, \phi) = \sum_{\alpha}U(\mu, \alpha)V(\alpha, \phi)\), and the choice of \(A_a = U\) and \(A_{{\Bar{a}}}= V\) forms a consistent choice of gauge between the pair \((C_a, C_{{\Bar{a}}})\).

Likewise, one can obtain \(Z_{b; \Bar{b}}\) and the SVD of \(Z_{b;\Bar{b}}\) results in \(A_b\) and \(A_{\Bar{b}}\). The same holds for
 \(Z_{\Bar{f};f}\), \(A_f\), and \(A_{\Bar{f}}\). To solve for \(G_{q}\) in the linear system \Cref{eqn: crossed linear system for G}, one simply contracts \(B_q\) with the pseudo-inverse of \(A_a, A_b, A_f\). Thus, one can use \Cref{eqn: crossed linear system for G} to solve for \(G_{q}\) for any \(l \not = 0, L\).

Our construction likewise gives rise to the equation of \(G_{q}\) for \(l = 0, L\) as special cases. For \(l = 0\), one can go through with the same calculation by setting \(C_f= S_f = A_f = 1\) in \Cref{eqn: uncrossed linear system for G} and \Cref{eqn: crossed linear system for G}.  For \(l = L\), one can sketch the linear system in \Cref{eqn: uncrossed linear system for G} for \(f = [d] - \{k\}\) by contraction with \(S_{f}\). The details for both cases are described in \Cref{alg:2}.

\paragraph{Algorithm summary}
The algorithm is summarized as \Cref{alg:2}, which demonstrates the whole procedure to carry out the functional hierarchical tensor sketching algorithm given a sample collection, including the details for the omitted special cases of \(l = 0, L\).

\begin{algorithm}[h!]
\caption{Functional hierarchical tensor interpolation.}
\label{alg:2}
\begin{algorithmic}[1]
\REQUIRE Function $g(x)$.
\REQUIRE Chosen function basis \(\{\psi_i\}_{i = 1}^{n}\).
\REQUIRE Random points 
\(\{x_{I_{k}^{(l)}}\}, \{x_{[d] - I_{k}^{(l)}}\}\)
generated from $\mathrm{Unif}([-1,1])$
and target internal ranks $\{r_{I_{k}^{(l)}}\}$ for each level index \(l\) and block index \(k\).
\FOR{each node \(q\) on the hierarchical tree}
    \STATE Set \(l\) as the level index of \(q\). Set \(k\) as the block index of \(q\).
    \IF{\(q\) is not leaf node}
        \STATE \((a, b) \gets (I_{2k-1}^{(l+1)}, I_{2k}^{(l+1)})\)
        \STATE Obtain \(Z_{a; \Bar{a}}, Z_{b; \Bar{b}}\) by \Cref{eqn: eq for Z_ag}
        \STATE Obtain \(A_{a}\) as the left factor of the best rank \(r_{a}\) factorization of \(Z_{a; \Bar{a}}\)
        \STATE Obtain \(A_{b}\) as the left factor of the best rank \(r_{b}\) factorization of \(Z_{b; \Bar{b}}\). 
    \ENDIF
    \IF{\(q\) is not root node}
        \STATE \(f \gets [d] -I_{k}^{(l)}\)
        \STATE Obtain \(Z_{\Bar{f}; f}\) by \Cref{eqn: eq for Z_ag}. 
        \STATE Obtain \(A_{f}\) as the right factor of the best rank \(r_{f}\) factorization of \(Z_{\Bar{f}; f}\).
    \ENDIF
    \IF{\(q\) is root node}
    \STATE Obtain \(
        B_{q}(\mu, \nu) = g(x_a^{\mu},x_b^{\nu}).
        \)
    \STATE Obtain \(G_{q}\) by solving the over-determined linear system \((A_{a} \otimes A_{b})G_{q} = B_{q}\).
    \ELSIF{\(q\) is leaf node}
    \STATE Obtain \(
        B_{q}(\mu, \zeta) = g(x_a^{\mu},x_f^{\zeta})
        \)
    \STATE Obtain \(\Psi_a\) as the function values of \(\{\psi_i\}_{i = 1}^{n}\) at \(\{x_a\}\).
    \STATE Obtain \(G_{q}\) by solving the over-determined linear system \((A_{f} \otimes \Psi_a)G_{q} = B_{q}\).
    \ELSE
    \STATE Obtain \(
        B_{q}(\mu, \nu, \zeta) = g(x_a^{\mu},x_b^{\nu},x_f^{\zeta})
        \)
    \STATE Obtain \(G_{q}\) by solving the over-determined linear system \((A_{a} \otimes A_{b} \otimes A_{f})G_{q} = B_{q}\).
    \ENDIF
\ENDFOR
\end{algorithmic}
\end{algorithm}


\section{Implementation details of FHT-based density estimation}\label{app: markov operator detail}

In this section, we describe the implementation details of computing an approximate joint density (Markov operator)
\(P(x,y,o) = \mathbb{P}\left[X = x, X' = y, O = o \right]\)
using FHT-based density estimation, where we note that in this section, the symbols \(O, o\) denote the random variable for the action variable to avoid confusion. We remark that the sketching mechanism in this section relies on sample moment estimation, whereas the sketching mechanism in \Cref{app: black-box interpolation detail} relies on function evaluation. The input to this routine is samples of \((X, X', O)\), and the output is the joint density \(P\).


We use the interlacing scheme to place the spatial variables on the leaf nodes.
We link \(o\) to the root node by an extra physical index such that the full tensor diagram represents the functional hierarchical tensor of \(g\).
A schematic diagram is shown in \Cref{fig:action_induced_Markov_operator_L_3}.
We denote the joint variable by \(z = (z_1, z_2, \ldots, z_{2d}, z_{2d+1}) = (x_1, y_1, \ldots, x_d, y_d, o)\).

As the construction for the variable hierarchical bipartition structure and the basis are fixed, one can apply the functional hierarchical tensor sketching algorithm to obtain an FHT-based characterization of the joint density function for \(Z = (X, X', O)\).  The derivation of hierarchical sketching with the introduction of an action variable does not significantly change the main routine of hierarchical sketching presented in Ref.~\cite{tang2024solvinga}. Mainly, one needs to modify the relevant sketch functions to take account of the action variable \(o\).


\paragraph{Equations}

In the FHT function class, the Markov operator $P$ admits a representation
\begin{equation}
    P(z)\equiv P(z_{1}, \ldots, z_{2d+1}) = \left<C, \, \bigotimes_{j=1}^{2d+1} \Vec{\Psi}_{j}(z_j) \right>,
\end{equation}
where \(\Vec{\Psi}_{j}\) is the collection of basis functions defined in \Cref{eqn: htn forward map} and $C$ is the coefficient tensor. The evaluation of \(g\) on the collected points \(\{z^{(i)}\}_{i = 1}^{N}\) reads \( P(z^{(i)}) = \left<C, \, \bigotimes_{j=1}^{2d+1} \Vec{\Psi}_{j}(z_j^{(i)}) \right> \).

We use \(G_{q}\) to denote a single tensor component, where \(q\) is the node of the \(k\)-th block at level \(l\).
For simplicity, we assume the general case $0< l < L$ so that \(q\) is neither the root nor the leaf node. We define \(a := I_{2k-1}^{(l+1)}, b := I_{2k}^{(l+1)}\) and \(f = [2d+1] - a \cup b\), where \(I_{k}^{(l)}\) follows notations introduced in \Cref{eqn: bipartition}. The structural low-rankness property of FHT implies that there exist \(r_a, r_b, r_f \in \mathbb{N}\), \(C_a \colon [n^{|a|}] \times [r_a] \to \R\), \(C_b \colon [n^{|b|}] \times [r_b] \to \R\), \(C_f \colon [r_f] \times [n^{|f|}] \to \R\) such that the following equations hold:
\begin{equation} \label{eqn: unsketched}
    C(i_a,i_b,i_f) = \sum_{\alpha, \beta, \theta} 
    C_{a}(i_{a}, \alpha)C_{b}(i_{b}, \beta)G_{q}(\alpha, \beta, \theta)C_{f}(\theta, i_{f}). 
\end{equation}

The hierarchical sketching algorithm essentially solves the over-determined linear system \Cref{eqn: unsketched} for \(G_{q}\) with the use of sketch functions. Let \(\tilde{r}_a, \tilde{r}_b, \tilde{r}_f\) be integers so that \(\tilde{r}_a > r_a, \tilde{r}_b > r_b, \tilde{r}_f > r_f\). Through the sketch functions \(S_a \colon [n^{|a|}] \times [\tilde{r}_a] \to \R, S_b \colon [n^{|b|}] \times [\tilde{r}_b] \to \R, S_f \colon [n^{|f|}] \times [\tilde{r}_f] \to \R\), one can contract the linear system in \Cref{eqn: unsketched} with \(S_a, S_b, S_f\) at the variables \(i_a, i_b, i_f\) respectively, which leads to the following sketched linear system:
\begin{equation} \label{eqn: sketched}
    B_q(\mu, \nu, \zeta) = \sum_{\alpha, \beta, \theta}A_{a}(\mu, \alpha)A_b(\nu, \beta)A_f(\zeta, \theta)G_{q}(\alpha, \beta, \theta),
\end{equation}
where \(A_{a}, A_b, A_f\) are respectively the contraction of \(C_a, C_b, C_f\) by \(S_a, S_b, S_f\) and \(B_q\) is the contraction of \(C\) by \(S_a \otimes S_b \otimes S_f\). The sketched-down linear system is thus:
\[
(A_{a} \otimes A_{b} \otimes A_{f})G_{q} = B_{q}.
\]

We define \(s_a \colon \R^{|a|} \times [\tilde{r}_a] \to \R\) by
\[
s_a(z_{a}, \mu) = \sum_{i_{j}, j \in a} S_a(i_{a}, \mu) \prod_{j \in a}\psi_{i_j}(z_j),
\]
and similar one defines \(s_b\) and \(s_f\) from \(S_b\) and \(S_f\). From the definition of \(s_a, s_b, s_f\), the term \(B_q\) satisfies
\[
B_q(\mu, \nu, \zeta) = \int_{[-1,1]^{2d+1}}  P(z_1, \ldots, z_{2d+1}) s_a(z_{a}, \mu)  s_b(z_{a}, \nu) s_f(z_{f}, \zeta)\, d \, z_1 \ldots z_{2d+1},
\]
and one has
\[
B_q(\mu, \nu, \zeta) = \mathbb{E}_{Z \sim P}\left[s_a(Z_a, \mu) s_b(Z_b, \nu) s_{f}(Z_f, \zeta)\right].
\]

As one only has access to samples of the probability density given by \(P\), one uses the approximation of \(P\) through empirical distribution, i.e. 
\(P \approx \frac{1}{N} \sum_{i=1}^N \delta_{z^{(i)}} \). Correspondingly, one obtains \(B_q\) through finite-sample moment estimation via the equation
\[ 
B_q(\mu, \nu, \zeta) = \frac{1}{N} \sum_{i=1}^N s_a(z_a^{(i)}, \mu) s_b(z_b^{(i)}, \nu) s_f(z_f^{(i)}, \zeta).
\]

To compute \(A_a,A_b,A_f\), the work in \cite{tang2024solvinga} shows that one similarly needs to compute quantities
\begin{equation} \label{eqn: Z_ag approximation}
    Z_{a; \Bar{a}}(\mu, \phi) \approx \frac{1}{N} \sum_{i = 1}^{N} s_a(z^{(i)}_{a}, \mu) s_{\Bar{a}}(z^{(i)}_{\Bar{a}}, \phi), 
\end{equation}
obtain \(Z_{a; \Bar{a}}(\mu, \phi) = \sum_\alpha U(\mu, \alpha) V(\alpha, \phi)\) from a singular value decomposition, and take \((A_a,A_{\Bar{a}}) = (U,V)\). The procedure also applies to the action-dependent case here. 
To solve for \(G_{q}\) in the linear system \Cref{eqn: sketched}, one simply contracts \(B_q\) with the pseudo-inverse of \(A_a, A_b, A_f\). Thus, one can use \Cref{eqn: sketched} to solve for \(G_{q}\) for any \(l \not = 0, L\).

Our construction likewise gives rise to the equation of \(G_{q}\) for \(l = 0, L\) as special cases. 
For \(l = 0\), one can go through with the same calculation by setting \(C_f= S_f = A_f = 1\) in \Cref{eqn: unsketched} and \Cref{eqn: sketched}.
For \(l=L\), one can sketch the linear system in \Cref{eqn: unsketched} for \(f = [2d+1] - \{k\}\) by contraction with \(S_{f}\). The details for both cases are described in \Cref{alg:3}.

\paragraph{Algorithm summary}
The algorithm is summarized as \Cref{alg:3}, which demonstrates the whole procedure to carry out the action-dependent operator-valued functional hierarchical tensor sketching algorithm given a sample collection.

\begin{algorithm}[h!]
\caption{Action-dependent operator-valued functional hierarchical tensor sketching.}
\label{alg:3}
\begin{algorithmic}[1]
\REQUIRE Sample \(\{z^{(i)}\}_{i = 1}^{N}\).
\REQUIRE Chosen function basis \(\{\psi_{i; j}\}_{i = 1}^{n}\) for each \(j \in [2d+1]\).
\REQUIRE Collection of sketch functions \(\{s_{I_{k}^{(l)}}\}, \{s_{[2d] - I_{k}^{(l)}}\}\) and target internal ranks $\{r_{I_{k}^{(l)}}\}$ for each level index \(l\) and block index \(k\).
\FOR{each node \(q\) on the hierarchical tree}
    \STATE Set \(l\) as the level index of \(q\). Set \(k\) as the block index of \(q\).
    \IF{\(q\) is not leaf node}
        \STATE \((a, b) \gets (I_{2k-1}^{(l+1)}, I_{2k}^{(l+1)})\)
        \STATE Obtain \(Z_{a; \Bar{a}}, Z_{b; \Bar{b}}\) by \Cref{eqn: Z_ag approximation}
        \STATE Obtain \(A_{a}\) as the left factor of the best rank \(r_{a}\) factorization of \(Z_{a; \Bar{a}}\)
        \STATE Obtain \(A_{b}\) as the left factor of the best rank \(r_{b}\) factorization of \(Z_{b; \Bar{b}}\). 
    \ENDIF
    \IF{\(q\) is not root node}
        \STATE \(f \gets [2d+1] -I_{k}^{(l)}\)
        \STATE Obtain \(Z_{\Bar{f}; f}\) by \Cref{eqn: Z_ag approximation}. 
        \STATE Obtain \(A_{f}\) as the right factor of the best rank \(r_{f}\) factorization of \(Z_{\Bar{f}; f}\).
    \ENDIF
    \IF{\(q\) is root node}
    \STATE Obtain \(
        B_{q}(\mu, \nu, m) = \frac{1}{N}\sum_{i}s_a(z^{(i)}_a, \mu) s_b(z^{(i)}_b, \nu) \psi_{2d+1;m}(z^{(i)}_{2d+1}).
        \)
    \STATE Obtain \(G_{q}\) by solving the over-determined linear system \((A_{a} \otimes A_{b})G_{q} = B_{q}\).
    \ELSIF{\(q\) is leaf node}
    \STATE Obtain \(
        B_{q}(j, \zeta) = \frac{1}{N}\sum_{i}\psi_{k;j}(z^{(i)}_k) s_{f}(z^{(i)}_f, \zeta).\)
    \STATE Obtain \(G_{q}\) by solving the over-determined linear system \((A_{f})G_{q} = B_{q}\).
    \ELSE
    \STATE Obtain \(
        B_{q}(\mu, \nu, \zeta) = \frac{1}{N}\sum_{i}s_a(z^{(i)}_a, \mu) s_b(z^{(i)}_b, \nu) s_f(z^{(i)}_f, \zeta).
        \)
    \STATE Obtain \(G_{q}\) by solving the over-determined linear system \((A_{a} \otimes A_{b} \otimes A_{f})G_{q} = B_{q}\).
    \ENDIF
\ENDFOR
\end{algorithmic}
\end{algorithm}

\end{document}